\documentclass{article}
 
\usepackage{amsmath,amssymb,amsfonts,amsthm,stmaryrd} 
\usepackage[colorlinks=true]{hyperref}
\usepackage[T1]{fontenc}
\usepackage[utf8]{inputenc} 
\usepackage{enumitem}
\usepackage{xcolor}
\usepackage{bbold} % indicatrice d'un ensemble
\usepackage{authblk}
\usepackage[a4paper]{geometry} 
\newtheorem{theorem}{Theorem}[section]

\newtheorem{lemma}[theorem]{Lemma}
\newtheorem{remark}[theorem]{Remark}

\newtheorem{proposition}[theorem]{Proposition}
\newtheorem{definition}[theorem]{Definition}  
\newtheorem{corollary}[theorem]{Corollary}

\newcommand{\vertiii}[1]{{\left\vert\kern-0.25ex\left\vert\kern-0.25ex\left\vert #1 
    \right\vert\kern-0.25ex\right\vert\kern-0.25ex\right\vert}}
\newcommand\numberthis{\addtocounter{equation}{1}\tag{\theequation}}
 
\hypersetup{
    colorlinks,
    linkcolor={blue!50!blue},
    citecolor={red}
}

\newtheorem{assumptionF}{Assumption} 
\newtheorem{assumptionO}{Assumption}

\title{Quasi-stationary distribution for the {L}angevin process in cylindrical domains, part I: existence, uniqueness and long-time convergence}
\author[1,2]{Tony Lelièvre\thanks{E-mail: tony.lelievre@enpc.fr}} 
\author[1,2]{Mouad Ramil\thanks{E-mail: mouad.ramil@enpc.fr}}
\author[1]{Julien Reygner\thanks{E-mail: julien.reygner@enpc.fr}}
\affil[1]{CERMICS, Ecole des Ponts, Marne-la-Vallée, France}
\affil[2]{MATHERIALS, Inria, Paris, France}
\date{\today} 

\begin{document}
\maketitle  
    
\begin{abstract}
Consider the Langevin process, described by a vector (position,momentum) in $\mathbb{R}^{d}\times\mathbb{R}^d$. Let $\mathcal O$ be a $\mathcal{C}^2$ open bounded and connected set of $\mathbb{R}^d$. We prove the compactness of the semigroup of the  Langevin process absorbed at the boundary of the domain $D:=\mathcal{O}\times\mathbb{R}^d$. We then obtain the existence of a unique quasi-stationary distribution (QSD) for the Langevin process on $D$.  We also provide a spectral interpretation of this QSD and obtain an exponential convergence of the Langevin process conditioned on non-absorption towards the QSD.   


\end{abstract}  
 
\section{Introduction} 
In statistical physics, the evolution of a molecular system at a given temperature is
typically modeled by the Langevin dynamics 
\begin{equation}\label{eq:Langevin_intro qsd}
  \left\{
    \begin{aligned}
        &\mathrm{d}q_t=M^{-1} p_t \mathrm{d}t , \\
        &\mathrm{d}p_t=F(q_t) \mathrm{d}t -\gamma M^{-1}  p_t
        \mathrm{d}t +\sqrt{2\gamma\beta^{-1}} \mathrm{d}B_t ,
    \end{aligned}
\right.  
\end{equation}
where $d=3N$ for a number $N$ of particles, $(q_t,p_t) \in \mathbb{R}^d \times \mathbb{R}^d$ denotes the set of positions and momenta of the
particles, $M \in \mathbb{R}^{d \times d}$ is the mass matrix,
$F:\mathbb{R}^d \to \mathbb{R}^d$ is the force acting on the particles, $\gamma >0$ is
the friction parameter, and $\beta^{-1}= k_B T$ with $k_B$
the Boltzmann constant and $T$
the temperature of the system. 
% Alternatively, the overdamped Langevin dynamics
% \begin{equation}\label{eq:ovLangevin_intro}
%   \mathrm{d}\overline{q}_t = F(\overline{q}_t)\mathrm{d} t + \sqrt{2\beta^{-1}}\mathrm{d} B_t,
% \end{equation}
% may also be employed. 
% Remarkably, both processes are related by the fact that when the force field $F$ is conservative, that is to say that there exists $V : \mathbb{R}^d \to \mathbb{R}$ such that $F = -\nabla V$, then the stationary distribution of $(\overline{q}_t)_{t \geq 0}$ writes
% \begin{equation}\label{eq:nuovLangevin}
%   \overline{\nu}(\mathrm{d}q) = \frac{1}{Z}\mathrm{e}^{-\beta V(q)}, \qquad Z = \int_{\mathbb{R}^d} \mathrm{e}^{-\beta V(q)} \mathrm{d}q,
% \end{equation}
% while the stationary distribution of $(q_t,p_t)_{t \geq 0}$ has the product structure
% \begin{equation}\label{eq:nuLangevin}
%   \nu(\mathrm{d}q\mathrm{d}p) = \overline{\nu}(\mathrm{d}q) \frac{\mathrm{e}^{-\frac{\beta\vert p\vert^2}{2}}}{(2\pi\beta^{-1})^{\frac{d}{2}}}\mathrm{d}p,
% \end{equation}
% the marginal in momentum of which is usually called the Maxwell distribution with inverse temperature $\beta$.

Such dynamics are used in particular to
compute thermodynamic and dynamic quantities, with numerous
applications in biology, chemistry and materials science. In many practical situations of interest, the system remains trapped for very long times
in subsets of the phase space, called metastable states, see for
example~\cite[Sections 6.3 and 6.4]{LelSto16}. This makes the simulation of these systems over the times of interest impossible. Typically, these states
are defined in terms of positions only, and are thus cylinders of the
form $D=\mathcal{O} \times \mathbb{R}^d$ for~\eqref{eq:Langevin_intro qsd}. In such a case, it is expected that the process reaches a local equilibrium distribution within the metastable state before leaving it. This distribution is called the quasi-stationary distribution (QSD). Proving the existence of this limiting behavior is in particular important to prove the consistency of accelerated dynamics algorithms, e.g. the parallel replica method, see for example~\cite{perez-uberuaga-voter-15}. It is also the building block to justify the use of jump Markov processes among the metastable states (kinetic Monte-Carlo or Markov state Models) to model the evolution over long timescales~\cite{kMC,MSM}.

While several works have already studied the properties of QSD for elliptic diffusion processes on a smooth bounded domain $\mathcal{O}$, to the best of our knowledge there are no available results for the Langevin dynamics~\eqref{eq:Langevin_intro qsd}, which is not elliptic but only hypoelliptic, and for which the natural domain $D = \mathcal{O} \times \mathbb{R}^d$ is not bounded, even if $\mathcal{O}$ is bounded. 
Building on several analytical results for the Langevin process~\eqref{eq:Langevin_intro qsd} obtained in~\cite{LelRamRey}, including a Gaussian upper-bound satisfied by the transition density of the Langevin process, we obtain the compactness of the semigroup of the Langevin process absorbed at the boundary of $D$. Applying the Krein-Rutman theorem we then obtain spectral properties on the infinitesimal generator of the Langevin process on $D$ with Dirichlet boundary conditions, and deduce the existence and uniqueness of a QSD $\mu$, as well as the fact that it describes the long-time behavior of the process conditioned on non absorption. 

Alternatively, a more probabilistic approach, based on general criteria developed by Champagnat and Villemonais~\cite{V}, is employed to obtain similar results in~\cite[Chapter 4]{RamPHD}. As we were finishing this work, we also became aware of the related work~\cite{GuiNectoux}, using different techniques based on Lyapounov functions. 
 
\medskip\textbf{Outline of the article.} In Section~\ref{sec: main results qsd}, we state the main results, which are then proven
in Section~\ref{section results qsd}.

\medskip\textbf{Notation.} Let us introduce here some  notation that will
be used in the following. We denote by $x=(q,p)$ generic elements in $\mathbb{R}^{2d}$, and  $\vert \cdot \vert$ the Euclidean norm both on $\mathbb{R}^d$ and on $\mathbb{R}^{2d}$. For a measurable subset $A$ of $\mathbb{R}^{2d}$, $\mathbb{R}_+^* \times \mathbb{R}^{2d}$ or $\mathbb{R}_+^* \times \mathbb{R}^{2d} \times \mathbb{R}^{2d}$, 
\begin{itemize} 
    \item $\vert A\vert$ is the Lebesgue measure of $A$,
    \item for $1 \le p\le \infty$, $\mathrm{L}^p(A)$ is the set of
      $\mathrm{L}^p$ scalar-valued functions on $A$ and $\Vert\cdot\Vert_{\mathrm{L}^p(A)}$ the associated norm,  
    \item $\mathcal{C}(A)$ (resp. $\mathcal{C}^b(A)$) is the set of
      scalar-valued continuous (resp. continuous and bounded)
      functions on $A$,  
    \item $\mathcal{C}^\infty(A)$ (resp.  $\mathcal{C}_c^\infty(A)$) is
      the set of scalar-valued $\mathcal{C}^\infty$ (resp. $\mathcal{C}^\infty$ with
      compact support) functions on $A$.  
\end{itemize}
We denote by $\|\cdot\|_\infty$ the sup norm on the Banach space $\mathcal{C}^b(A)$. For $T$ a linear bounded operator on $\mathcal{C}^b(A)$, we denote its operator norm by:
$$\vertiii{T}_{\mathcal{C}^b(A)}:=\sup_{f\in \mathcal{C}^b(A),\Vert f\Vert_\infty\leq1}\Vert Tf\Vert_\infty.$$

\section{Main results}\label{sec: main results qsd}

This section presents the main results we obtained. 

As a motivation, we first recall in Section~\ref{subsection elliptic} what is known about the QSD of the overdamped Langevin process. 

In order to prepare the presentation of our main results, we state various analytical properties of the Langevin process and the related kinetic Fokker-Planck equation in Section~\ref{subsection state of the art}. The proofs of these auxiliary results are detailed in~\cite{LelRamRey}. 

Our main results, concerning the degenerate case of the Langevin process, are presented in Section~\ref{Main results_QSD}. We first state the compactness of the semigroup of the Langevin process absorbed at the boundary of $D$. The existence of a unique QSD of the Langevin process in $D$ is then obtained. Besides, this QSD is shown to be the unique solution of an eigenvalue problem related to the infinitesimal generator of the process $(q_t,p_t)_{t \ge 0}$ with absorbing boundary conditions. Finally, this QSD attracts all probability measures on $D$, at an exponential rate. 
 
 \medskip
Let us conclude this introduction by recalling the definitions of the quasi-stationary and quasi-limiting distributions, which are the central notions of this work, in a general setting. We refer to~\cite{Collet,VilMel} for a complete introduction. 
Let $E$ be a Polish space endowed with its Borel $\sigma$-algebra $\mathcal{B}(E)$, and let $(X_t)_{t \geq 0}$ be a time-homogeneous, strong Markov process in $E$ with continuous sample-paths. For any $x \in E$, we denote by $\mathbb{P}_x$ the probability measure under which $X_0=x$ almost surely, and for any probability measure $\theta$ on $E$, we define
\begin{equation*}
  \mathbb{P}_\theta(\cdot) := \int_E \mathbb{P}_x(\cdot)\theta(\mathrm{d}x).
\end{equation*}
Let $D$ be an open subset of $E$ and $\tau_\partial$ be the stopping time defined by
\begin{equation*}
  \tau_\partial := \inf\{t>0: X_t \not\in D\}.
\end{equation*}

\begin{definition}[QSD]\label{def QSD 1} A probability measure $\mu$ on $D$ is said to be a QSD on $D$ of the process $(X_t)_{t\geq0}$, if for all $A\in\mathcal{B}(D) := \{A\cap D,A\in\mathcal{B}(E)\}$, for all $t\geq0$,
\begin{equation}\label{eq:defqsd}
  \mathbb{P}_\mu(X_t\in A, \tau_\partial>t)=\mu(A)\mathbb{P}_\mu(\tau_\partial>t).
\end{equation}
\end{definition}
When $\mathbb{P}_\mu(\tau_\partial>t)>0$, the identity~\eqref{eq:defqsd} equivalently writes $\mathbb{P}_\mu(X_t\in A|\tau_\partial>t)=\mu(A)$.
 
A closely related notion is that of Quasi-Limiting Distribution (QLD), which is a probability measure $\mu$ on $D$ such that there exists a probability measure $\theta$ on $D$ for which
\begin{equation}\label{eq:qld}
  \forall A \in \mathcal{B}(D), \qquad \mu(A) = \lim_{t \to \infty} \mathbb{P}_\theta(X_t \in A | \tau_\partial > t).    
\end{equation}
A QLD is necessarily a QSD (and the converse is obvious), and we say that $\mu$ "attracts" $\theta$ when~\eqref{eq:qld} holds. When a QSD attracts all Dirac masses on $D$, it is called a Yaglom limit.

Last, when $X_0$ is initially distributed according to a QSD on $D$, the exit event from $D$ of the process $(X_t)_{t\geq0}$ satisfies the following properties, see~\cite[Theorems 2.2 and 2.6]{Collet}.
\begin{proposition}[Exit event]\label{exit event general case}
Let $\mu$ be a QSD on $D$ of the process $(X_t)_{t\geq0}$, then there exists $\lambda_0\geq0$ such that
\begin{enumerate}
    \item $\tau_\partial$ follows the exponential law of parameter $\lambda_0$, that is to say $\mathbb{P}_\mu(\tau_\partial>t)=\mathrm{e}^{-\lambda_0 t}$ for all $t \geq 0$,
    \item if $\lambda_0>0$, $X_{\tau_\partial}$ is independent of $\tau_\partial$. 
\end{enumerate} 
\end{proposition} 
In the former statement, the case $\lambda_0 = 0$ means that $\tau_\partial=\infty$, $\mathbb{P}_\mu$-almost surely.

\subsection{Elliptic case and the overdamped Langevin process}\label{subsection elliptic}
 Quasi-stationary distributions on smooth bounded domains for elliptic diffusion processes, have been widely studied in the literature. We refer for example to~\cite{GQZ,LebLelPer,V2,V3}. Let us recall here some of their important results. 

Let $\beta>0$ and $F:\mathbb{R}^d\mapsto\mathbb{R}^d$ satisfying the following assumption. 
\begin{assumptionF}\label{hyp F1 qsd}
$F\in\mathcal{C}^\infty(\mathbb{R}^{d},\mathbb{R}^{d})$. 
\end{assumptionF}

Let $(\Omega,\mathcal{F},(\mathcal{F}_t)_{t\geq0},\mathbb{P})$ be a filtered probability space and $(B_t)_{t\geq0}$ a $d$-dimensional $(\mathcal{F}_t)_{t\geq0}$-Brownian motion. Consider the overdamped Langevin process defined by
\begin{equation}\label{eq:ovLangevin_intro}
  \mathrm{d}\overline{q}_t = F(\overline{q}_t)\mathrm{d} t + \sqrt{2\beta^{-1}}\mathrm{d} B_t.
\end{equation}

Under Assumption~\ref{hyp F1 qsd}, the vector field $F$ is locally Lipschitz continuous and therefore the stochastic differential equation~\eqref{eq:ovLangevin_intro} possesses a unique strong solution $(\overline{q}_t)_{0 \leq t < \overline{\tau}_\infty}$ defined up to some explosion time $\overline{\tau}_\infty \in (0,+\infty]$. The overdamped Langevin process admits the following infinitesimal generator: 
\begin{equation}\label{L}
    \overline{\mathcal{L}}=F\cdot\nabla+\beta^{-1}\Delta, 
\end{equation}
with formal adjoint $\overline{\mathcal{L}}^*$ in $\mathrm{L}^2(\mathrm{d}x)$ given by:  
\begin{equation*}\label{generateur adjoint overdamped}
    \overline{\mathcal{L}}^*=-\mathrm{div}(F\cdot )+\beta^{-1}\Delta. 
\end{equation*}
Let $\mathcal{O}$ be an open set of $\mathbb{R}^d$ satisfying the following assumption.
\begin{assumptionO}\label{hyp O qsd}
$\mathcal{O}$ is an open $\mathcal{C}^2$ bounded connected set of $\mathbb{R}^d$. 
\end{assumptionO}
Let $\overline{\tau}_{\partial}:=\inf \{t>0: \overline{q}_t\notin \mathcal{O}\}$ be the first exit time from $\mathcal{O}$ of the process $(\overline{q}_t)_{0 \leq t < \overline{\tau}_\infty}$. Under Assumption~\ref{hyp O qsd}, the vector field $F$ is Lipschitz continuous on $\mathcal{O}$ and therefore $\overline{\tau}_\partial \leq \overline{\tau}_\infty$. 

It has been shown in~\cite{V3,GQZ,LebLelPer,KnobPart} that the overdamped Langevin process admits a unique QSD on $\mathcal{O}$, which moreover satisfies the following properties.
\begin{theorem}[QSD of the overdamped Langevin process]\label{qsd overdamped} Under Assumptions~\ref{hyp F1 qsd} and~\ref{hyp O qsd}, there exists a unique QSD $\overline{\mu}$ on $\mathcal{O}$ of the process $(\overline{q}_t)_{t\geq0}$. Furthermore, 
\begin{enumerate}[label=(\roman*),ref=\roman*]
    \item there exists $\overline{\psi}\in\mathcal{C}^2(\mathcal{O})\cap\mathcal{C}^b(\overline{\mathcal{O}})$ such that $\overline{\mu}(\mathrm{d}q)=\overline{\psi}(q)\mathrm{d}q$, where $\mathrm{d}q$ is the Lebesgue measure on $\mathbb{R}^d$,
    \item $\mathrm{Span}(\overline{\psi})$ is the eigenspace associated with the smallest eigenvalue $\overline{\lambda}$ of the operator $-\overline{\mathcal{L}}^*$ with homogeneous Dirichlet boundary conditions on $\partial\mathcal{O}$, 
    \item there exist $C>0$ and $\alpha>0$ such that for all probability measures $\theta$ on $\mathcal{O}$, for all $t\geq0$, $$\big\Vert \mathbb{P}_\theta(\overline{q}_t\in\cdot\vert\overline{\tau}_\partial>t)-\overline{\mu}(\cdot)\big\Vert_{TV}\leq C\mathrm{e}^{-\alpha t},$$ where $\Vert\cdot\Vert_{TV}$ is the total-variation norm on the space of bounded signed measures on $\mathbb{R}^d$.
\end{enumerate}  
\end{theorem}
 
Multiple approaches are used in the literature to obtain the properties above. In the conservative case $F=-\nabla V$, under suitable assumptions on $V$ we have $\overline{\tau}_\infty = \infty$ and the process $(\overline{q}_t)_{t\geq0}$ is reversible with respect to the measure $\mathrm{e}^{-\beta V(q)}\mathrm{d}q$. As a consequence, $\overline{\mathcal{L}}$ is symmetric with respect to the canonical scalar product on $\mathrm{L}^2(\mathrm{e}^{-\beta V(q)}\mathrm{d}q)$ and since the inverse of the operator $\overline{\mathcal{L}}$ with homogeneous Dirichlet boundary condition on $\partial\mathcal{O}$ is compact from $\mathrm{L}^2(\mathrm{e}^{-\beta V(q)}\mathrm{d}q)$ to $\mathrm{L}^2(\mathrm{e}^{-\beta V(q)}\mathrm{d}q)$, one can obtain a discrete spectral decomposition of $\overline{\mathcal{L}}$ with this boundary condition. This then yields the theorem above, see
~\cite{LebLelPer}. 

In the general case when $F$ is non conservative, the process $(\overline{q}_t)_{t\geq0}$ is not necessarily reversible but a spectral approach can still be used. In~\cite{GQZ} the authors prove the compactness of the semigroup $(\overline{P}^\mathcal{O}_t)_{t \geq 0}$ defined on the Banach space  
\begin{equation*}
  \{f\in\mathcal{C}^b(\mathcal{O}) : \forall q\in\mathcal{O}, f(q)=\mathrm{d}_\partial(q) g(q)\quad\text{s.t. }g\text{ is uniformly continuous on }\mathcal{O}\},
\end{equation*}
where $\mathrm{d}_\partial$ is the Euclidean distance to the boundary $\partial\mathcal{O}$, by
\begin{equation*}
  \overline{P}^\mathcal{O}_tf:x\in D\mapsto \mathbb{E}_q\left[f(\overline{q}_t)\mathbb{1}_{\overline{\tau}_\partial>t}\right],
\end{equation*} 
using sharp estimates of the Green function of $\overline{\mathcal{L}}$ shown in~\cite{GrutWid}. Then, applying Krein-Rutman theorem to the operator $\overline{P}^\mathcal{O}_t$, the authors manage to deduce Theorem~\ref{qsd overdamped}.

Last, a more probabilistic approach is developed in~\cite{V3} where the authors prove that the semigroup $(\overline{P}^\mathcal{O}_t)_{t \geq 0}$ satisfies a gradient estimate, irreducibility conditions and a controlled probability of absorption near the boundary $\partial\mathcal{O}$ which also yields Theorem~\ref{qsd overdamped}.

\subsection{Analytical properties of the Langevin process} \label{subsection state of the art}
In this section we recall some results from~\cite{LelRamRey} that will be used henceforth. Let $\gamma\in\mathbb{R}$, $\sigma>0$. Under Assumption~\ref{hyp F1 qsd}, the stochastic differential equation 
\begin{equation}\label{Langevin qsd}
\left\{
    \begin{aligned}
        \mathrm{d}q_t&=p_t \mathrm{d}t , \\
        \mathrm{d}p_t&=F(q_t) \mathrm{d}t -\gamma  p_t \mathrm{d}t+\sigma\mathrm{d}B_t,
    \end{aligned}
  \right.
\end{equation}
possesses a unique strong solution $(X_t=(q_t,p_t))_{0 \leq t < \tau_\infty}$, defined up to some explosion time $\tau_\infty \in (0,+\infty]$.
Notice that, compared to~\eqref{eq:Langevin_intro qsd}, we
consider here and
  henceforth the mass to be identity without loss of generality (see
the change of variables in~\cite[Equation
(3.117)]{LelRouSto10}), so that momentum is identified with velocity. Besides, we consider the general case $\gamma\in\mathbb{R}$ and
$\sigma>0$ not necessarily related to $\gamma$.
 
The infinitesimal generator of the Langevin process is the operator $\mathcal{L}$, defined for $(q,p)\in\mathbb{R}^d\times\mathbb{R}^d$ by:
\begin{equation}\label{generateur Langevin qsd}
    \mathcal{L} =  p\cdot\nabla_q+F(q)\cdot\nabla_p -\gamma  p\cdot\nabla_p+ \frac{\sigma^2}{2}\Delta_p, 
\end{equation}
with formal adjoint $\mathcal{L}^*$ in $\mathrm{L}^2(\mathrm{d}x)$ given by:
\begin{equation}\label{generateur adjoint qsd}
    \mathcal{L}^*=-p\cdot\nabla_q -F(q)\cdot\nabla_p+\gamma \mathrm{div}_p(p  \cdot )+\frac{\sigma^2}{2}\Delta_p. 
\end{equation}

Consider now the following strenghtening of Assumption~\ref{hyp F1 qsd}.
\begin{assumptionF}\label{hyp F2 qsd}
$F\in\mathcal{C}^\infty(\mathbb{R}^{d},\mathbb{R}^{d})$ and $F$ is bounded and globally Lipschitz continuous on $\mathbb{R}^d$. 
\end{assumptionF} 
Under Assumption~\ref{hyp F2 qsd}, $\tau_\infty=\infty$ almost surely and the Langevin process~\eqref{Langevin qsd} admits a smooth transition density $\mathrm{p}_t(x,y)$~\cite[Corollary 7.2]{RB}, which is positive~\cite[Corollary 3.3]{PositiveDensity}. In addition, this density admits an explicit Gaussian upper-bound, see~\cite[Theorem 2.19]{LelRamRey}.  

\begin{theorem}[Gaussian upper-bound]\label{borne densite thm qsd}
Under Assumption~\ref{hyp F2 qsd}, the transition density $\mathrm{p}_t(x,y)$ of the Langevin process $(X_t)_{t\geq0}$ satisfying~\eqref{Langevin qsd} is such that for all $\alpha\in (0,1)$, there exists $c_\alpha>0$, depending only on $\alpha$, such that for all $T>0$ and $t\in(0,T]$, for all $x,y\in\mathbb{R}^{2d}$,
\begin{equation} \label{borne densité qsd}
    \mathrm{p}_t(x,y)\leq C_{t,T} \widehat{\mathrm{p}}^{(\alpha)}_{t}(x,y),
\end{equation}
where $C_{t,T}:=\frac{1}{\alpha^d} \sum_{j=0}^\infty \frac{\left( \Vert F\Vert_\infty c_\alpha(1+\sqrt{\gamma_- T}) \sqrt{\pi t}  \right)^j}{\sigma^j \Gamma\left(\frac{j+1}{2}\right)}$, $\gamma_-=\max(-\gamma,0)$ is the negative part of $\gamma\in\mathbb{R}$,
$\Gamma$ is the Gamma function and
$\widehat{\mathrm{p}}^{(\alpha)}_{t}(x,y)$ is the transition density
of the Gaussian process
 $(\widehat{q}^{(\alpha)}_t, \widehat{p}^{(\alpha)}_t)_{t \geq 0}$ defined by
\begin{equation}\label{eq:processus alpha}
  \left\{\begin{aligned}
    \mathrm{d}\widehat{q}^{(\alpha)}_t &= \widehat{p}^{(\alpha)}_t \mathrm{d}t,\\
    \mathrm{d}\widehat{p}^{(\alpha)}_t &= -\gamma\widehat{p}^{(\alpha)}_t \mathrm{d}t + \frac{\sigma}{\sqrt{\alpha}}\mathrm{d}B_t.
  \end{aligned}\right.
\end{equation}  
\end{theorem}  

We now let $\mathcal{O} \subset \mathbb{R}^d$ satisfy Assumption~\ref{hyp O qsd} and consider the following domain of $\mathbb{R}^{2d}$,  $$D:=\mathcal{O}\times\mathbb{R}^d,$$ where the first coordinate (position) is constrained to remain on the bounded open set $\mathcal{O}$ and the second one (velocity) remains free. This
is the natural phase space domain of the Langevin process absorbed when
leaving $\mathcal{O}$. 

For $q\in\partial\mathcal{O}$, let
$n(q)\in\mathbb{R}^d$ be the unitary outward normal vector to
$\mathcal{O}$ at $q \in \partial\mathcal{O}$. We introduce the following partition of $\partial D$: 
$$ \Gamma^+=\{ (q,p)\in\partial\mathcal{O}\times\mathbb{R}^d : p\cdot n(q)>0 \} ,$$
$$ \Gamma^-=\{ (q,p)\in\partial\mathcal{O}\times\mathbb{R}^d : p\cdot n(q)<0 \} ,$$
$$ \Gamma^0=\{ (q,p)\in\partial\mathcal{O}\times\mathbb{R}^d : p\cdot n(q)=0 \}.$$
 
Let $\tau_{\partial}$ be the first exit time from $D$ of the Langevin process $(X_t)_{t\geq0}$ in~\eqref{Langevin qsd}, i.e.  $$\tau_{\partial}=\inf \{t>0: X_t\notin D\} .$$ 
Under Assumptions~\ref{hyp F1 qsd} and~\ref{hyp O qsd}, $F$ is Lipschitz continuous on $\mathcal{O}$ and therefore $\tau_\partial \leq \tau_\infty$.
\begin{remark}\label{solution dans D}
  Friedman's uniqueness result~\cite[Theorem 5.2.1.]{F}
ensures that the trajectories $(X_t)_{0 \le t \le \tau_\partial}$ do not depend on the values of $F$ outside of $\mathcal{O}$. Therefore, whenever we are interested in quantities which only depend on the absorbed Langevin process, there is no loss of generality in modifying $F$ outside of $\mathcal{O}$ so that it satisfies Assumption~\ref{hyp F2 qsd}.
\end{remark}  
 The Langevin process absorbed outside of the domain $D$ has been thoroughly studied in~\cite{LelRamRey}. Some of the results associated to its transition density are reminded below and can be found in~\cite[Theorem 2.20, Corollary 2.21]{LelRamRey}.
\begin{theorem}[Transition density of the absorbed Langevin process]\label{thm density intro qsd}
    Under Assumptions~\ref{hyp F1 qsd} and~\ref{hyp O qsd}, there exists a function $$(t,x,y)\mapsto \mathrm{p}_t^D(x,y) \in\mathcal{C}^\infty(\mathbb{R}_+^{*}\times D\times D)\cap\mathcal{C}(\mathbb{R}_+^*\times\overline{D}\times\overline{D})$$ 
    which satisfies for all $t>0$,  
\begin{itemize}
    \item $\mathrm{p}^D_t(x,y)>0$ for all $x\notin \Gamma^+\cup\Gamma^0$ and $y\notin\Gamma^-\cup\Gamma^0$, 
    \item $\mathrm{p}^D_t(x,y)=0$ if $x\in \Gamma^+\cup\Gamma^0$ or if $y\in\Gamma^-\cup\Gamma^0$,
\end{itemize}
    and is such that for all $t>0$, $x\in \overline{D}$ and $A\in\mathcal{B}(D)$, $$\mathbb{P}_x(X_t\in A,\tau_\partial>t)=\int_A \mathrm{p}_t^D(x,y) \mathrm{d}y.$$ Moreover, for $f\in\mathcal{C}
   ^b(\overline{D})$, the functions $u,v$ defined by:   $$\forall t>0,\quad \forall x\in D,\qquad u(t,x):=\int_D \mathrm{p}_t^D(x,y)f(y)\mathrm{d}y,\qquad v(t,x):=\int_D \mathrm{p}_t^D(y,x)f(y)\mathrm{d}y,$$ are in $\mathcal{C}^\infty(\mathbb{R}_+^*\times D)$ and satisfy:
   $$\forall t>0,\quad \forall x\in D,\qquad\partial_tu(t,x)=\mathcal{L}u(t,x),\qquad\partial_tv(t,x)=\mathcal{L}^*v(t,x).$$ 
   Finally, $\mathrm{p}^D_t(x,y)$ also satisfies the Gaussian upper-bound~\eqref{borne densité qsd} where $\|F\|_\infty$ is replaced with $\|F\|_{\mathrm{L}^\infty(D)}$ in $C_{t,T}$.
\end{theorem}
 
We conclude this subsection with a time-reversibility result from~\cite[Section 6.1]{LelRamRey} linking the transition densities of the Langevin process~\eqref{Langevin qsd} and of the process called "adjoint" Langevin process $(\Tilde{X}_t=(\Tilde{q}_t,\Tilde{p}_t))_{t\geq0}$ with infinitesimal generator $\Tilde{\mathcal{L}}:=\mathcal{L}^*-d\gamma$, and satisfying the following SDE:  
\begin{equation}\label{Langevin adjoint qsd}
  \left\{
    \begin{array}{ll}
        \mathrm{d}\Tilde{q}_t=-\Tilde{p}_t \mathrm{d}t , \\
        \mathrm{d}\Tilde{p}_t=-F(\Tilde{q}_t) \mathrm{d}t+\gamma \Tilde{p}_t \mathrm{d}t+\sigma \mathrm{d}B_t.
    \end{array}
\right.   
\end{equation} 
Let $\Tilde{\tau}_{\partial}$ be the first exit time from $D$ of $\Tilde{X}_t$, i.e. $\tilde{\tau}_{\partial}:=\inf \{t>0: \Tilde{X}_t\notin D\}$. The transition kernel $\mathbb{P}_x(\Tilde{X}_t\in\cdot,\Tilde{\tau}_\partial>t)$ admits a transition density $\Tilde{\mathrm{p}}_t^D(x,y)$ which satisfies the following equality, see~\cite[Theorem 6.2]{LelRamRey}.
\begin{theorem}[Time-reversibility]\label{duality thm qsd} Under Assumptions~\ref{hyp F1 qsd} and~\ref{hyp O qsd}, 
\begin{equation} 
    \forall t>0, \quad \forall x,y\in D,\qquad\mathrm{p}_t^D(x,y)=\mathrm{e}^{d \gamma t} \Tilde{\mathrm{p}}_t^D(y,x).
\end{equation} 
\end{theorem} 
 
\subsection{Compactness and QSD of the Langevin process} \label{Main results_QSD}
In this section we state the main results proven in this work. Let us emphasize the fact that the results stated in the present section hold for any $\gamma \in \mathbb{R}$, $\sigma>0$ and $F$ satisfying Assumption~\ref{hyp F1 qsd}. The first result states the compactness of the semigroup $(P^D_t)_{t\geq0}$, defined below, based on the Gaussian estimate from Theorem~\ref{thm density intro qsd}, which implies that $\mathrm{p}^D_t\in\mathrm{L}^\infty(D\times D)\cap\mathrm{L}^1(D\times D)$ for all $t>0$ (see Lemma~\ref{hilbert schmidt lemma}).
\begin{theorem} [Semigroup of the absorbed Langevin process]\label{compactness} Let Assumptions~\ref{hyp O qsd} and~\ref{hyp F1 qsd} hold. For any $t \geq0$, $p \in [1,+\infty]$ and $f \in \mathrm{L}^p(D)$, the quantity
\begin{equation}\label{def semigroupe}
  P^D_t f : x \in \overline{D} \mapsto \mathbb{E}_x\left[\mathbb{1}_{\tau_\partial > t} f(X_t)\right]
\end{equation} 
is well-defined. Besides, let $p,q\in [1,+\infty]$.
\begin{enumerate}[label={\rm(\roman*)},ref=\roman*]
  \item\label{it:compacite-qsd:1} The family of operators $(P^D_t)_{t \geq 0}$ is a semigroup on $\mathrm{L}^p(D)$ and on $\mathcal{C}^{b}(\overline{D})$. 
  \item\label{it:compacite-qsd:2} For any $t>0$, the operator $P^D_t$ maps $L^p(D)$ into $L^q(D)$ and into $\mathcal{C}^b(\overline{D})$ continuously. 
  \item\label{it:compacite-qsd:3} For any $t>0$, the operator $P^D_t$ is compact from $\mathrm{L}^p(D)$ to $\mathrm{L}^p(D)$, and from $\mathcal{C}^{b}(\overline{D})$ to $\mathcal{C}^{b}(\overline{D})$. 
\end{enumerate}
\end{theorem}
\begin{remark}
  As an alternative to our probabilistic approach, we expect that a similar statement might also be deduced from the subelliptic estimates on the kinetic Fokker-Planck operator recently obtained by Nier~\cite{nier-18} in a very general framework (both on the geometry of the underlying phase space and on the boundary conditions). 
\end{remark}
Similarly, one can define the family of operators $(\Tilde{P}^D_t)_{t\geq0}$ associated to the transition density $\Tilde{\mathrm{p}}_t^D$, given for $f\in\mathrm{L}^\infty(D)$ by 
\begin{equation}\label{def semigroupe adjoint}
  \Tilde{P}^D_t f : x \in \overline{D} \mapsto \mathbb{E}_x\left[\mathbb{1}_{\Tilde{\tau}_\partial > t} f(\Tilde{X}_t)\right].
\end{equation} 
\begin{remark}\label{semigroupe adjoint}
  It follows from Theorem~\ref{duality thm qsd} and Lemma~\ref{hilbert schmidt lemma} that $\tilde{\mathrm{p}}^D_t\in\mathrm{L}^\infty(D\times D)\cap\mathrm{L}^1(D\times D)$ for all $t>0$. Therefore, following the proof of Theorem~\ref{compactness} in Section~\ref{sec:compact}, one can also obtain that the family $(\Tilde{P}^D_t)_{t\geq0}$ defined in~\eqref{def semigroupe adjoint} satisfies the properties detailed in Theorem~\ref{compactness}. In particular, $\Tilde{P}^D_t:\mathcal{C}^b(\overline{D})\to\mathcal{C}^b(\overline{D})$ and $\Tilde{P}^D_t:\mathrm{L}^p(D)\to\mathrm{L}^p(D)$ for $p \in [1,+\infty]$ are compact.
\end{remark}

Our second result focuses on the spectral radii of $(P^D_t)_{t\geq0}$ and $(\Tilde{P}^D_t)_{t\geq0}$. Let us recall here the definition of the spectral radius of a bounded operator $T$ on the Banach space $\mathcal{C}^b(\overline{D})$, which can be found in~\cite[p.~192]{RS}.

\begin{definition}[Spectrum and spectral radius]\label{def spectral radius}
Let $T$ be a bounded real operator on $\mathcal{C}^b(\overline{D})$ and $I$ the identity operator. Let us call $\sigma(T)$ the spectrum of $T$ which is defined by:
$$\sigma(T):=\{\lambda\in\mathbb{C}:T-\lambda I\text{ does not have an inverse that is a bounded linear operator}\}.$$
The spectral radius $r(T)$ of $T$ is then defined as:
$$r(T):=\sup_{\lambda\in\sigma(T)}\vert\lambda\vert .$$ 
\end{definition}
We obtain the following result on the operators $P^D_t$ and $\Tilde{P}^D_t$ (defined in~\eqref{def semigroupe} and~\eqref{def semigroupe adjoint}) and their spectral radius.
\begin{theorem}[Spectral properties of $P^D_t$ and $\Tilde{P}^D_t$]\label{properties phi}
Under Assumptions~\ref{hyp F1 qsd} and~\ref{hyp O qsd}, there exists $\lambda_0>0$ such that for all $t\geq 0$,
$$r(P^D_t)=\mathrm{e}^{-\lambda_0t},\qquad r(\Tilde{P}^D_t)=\mathrm{e}^{-(\lambda_0+d\gamma)t}.$$
Besides, there exist unique functions $\phi,\psi\in\mathcal{C}^b(\overline{D})$, up to a multiplicative constant, such that for all $t \geq 0$,
$$P^D_t\phi=\mathrm{e}^{-\lambda_0 t}\phi\quad\text{and}\quad \Tilde{P}^D_t\psi=\mathrm{e}^{-(\lambda_0+d\gamma) t}\psi . $$
Last, $\phi,\psi\in\mathrm{L}^1(D)\cap\mathcal{C}^\infty(D)$ and 
\begin{itemize}
    \item $\phi>0$ on $D\cup\Gamma^-$, $\phi=0$ on $\Gamma^+\cup\Gamma^0$ and $\mathcal{L}\phi=-\lambda_0 \phi$ on $D$,
    \item $\psi>0$ on $D\cup\Gamma^+$, $\psi=0$ on $\Gamma^-\cup\Gamma^0$ and $\mathcal{L}^*\psi=-\lambda_0 \psi$ on $D$.
\end{itemize} 
\end{theorem}

In the following, we choose $\phi$ and $\psi$ such that $\int_D\phi(x)\mathrm{d}x=\int_D\psi(x)\mathrm{d}x=1$. The proof of this theorem relies on the application of the Krein-Rutman theorem~\cite[p. 313]{SW} to the compact operators $P^D_t$ and $\Tilde{P}^D_t$ on the Banach space $\mathcal{C}^b(\overline{D})$. We are able to deduce from this result the existence of a unique QSD on $D$ for the processes $(X_t)_{t\geq0}$ and $(\Tilde{X}_t)_{t\geq0}$. 
\begin{theorem}[Existence and uniqueness of a QSD] \label{density qsd thm}Let Assumptions~\ref{hyp F1 qsd} and~\ref{hyp O qsd} hold. Let $\mu$ and $\Tilde{\mu}$ be the following probability measures on $D$: 
\begin{equation}\label{densite qsd}
    \forall A\in\mathcal{B}(D),\qquad\mu(A)=\int_A \psi(x)\mathrm{d}x,\quad\Tilde{\mu}(A)=\int_A \phi(x)\mathrm{d}x.
\end{equation}
Then $\mu$ (resp. $\Tilde{\mu}$) is the unique QSD on $D$ of $(X_t)_{t\geq0}$ (resp. $(\Tilde{X}_t)_{t\geq0}$) and for all $t \geq 0$,
\begin{equation}\label{law exit time}
    \mathbb{P}_\mu(\tau_\partial>t)=\mathrm{e}^{-\lambda_0 t},\quad\mathbb{P}_{\Tilde{\mu}}(\Tilde{\tau}_\partial>t)=\mathrm{e}^{-(\lambda_0+d \gamma) t}.
\end{equation} 
\end{theorem} 
Moreover, one can characterize the law of $X_{\tau_\partial}$ (resp. $\tilde{X}_{\tilde{\tau}_\partial}$) when $X_0\sim\mu$ (resp. $\tilde{X}_0\sim\tilde{\mu}$) using the following theorem.
\begin{theorem}\label{prop: exit event qsd}
For all $t\geq0$, $f\in\mathrm{L}^\infty(\partial D)$,
\begin{align*}\label{eq exit point law}
    \int_D\psi(q,p)\mathbb{E}_{(q,p)}\left[f(q_{\tau_\partial},p_{\tau_\partial})\mathbb{1}_{\tau_\partial\leq t}\right]\mathrm{d}q\mathrm{d}p&=\frac{1-\mathrm{e}^{-\lambda_0t}}{\lambda_0}\int_{\partial D}\psi(\overline{q},p)f(\overline{q},p)\left\vert p\cdot n(\overline{q})\right\vert\sigma_{\partial \mathcal{O}}(\mathrm{d}\overline{q})\mathrm{d}p,\\
    \int_D\phi(q,p)\mathbb{E}_{(q,p)}\left[f(\tilde{q}_{\tilde{\tau}_\partial},\tilde{p}_{\tilde{\tau}_\partial})\mathbb{1}_{\tilde{\tau}_\partial\leq t}\right]\mathrm{d}q\mathrm{d}p&=\frac{1-\mathrm{e}^{-(\lambda_0+d\gamma)t}}{\lambda_0+d\gamma}\int_{\partial D}\phi(\overline{q},p)f(\overline{q},p)\left\vert p\cdot n(\overline{q})\right\vert\sigma_{\partial \mathcal{O}}(\mathrm{d}\overline{q})\mathrm{d}p,
\end{align*}
where $\sigma_{\partial \mathcal{O}}$ is the surface measure on $\partial \mathcal{O}$.
\end{theorem}
In fact, taking $t\rightarrow\infty$, one obtain the following description of the first exit event starting from the QSD on $D$.
\begin{corollary}[First exit point starting from the QSD]
If $X_{0}\sim\mu$, 
$$X_{\tau_\partial}\sim\frac{1}{\lambda_0}\left\vert p\cdot n(\overline{q})\right\vert\psi(\overline{q},p)\sigma_{\partial \mathcal{O}}(\mathrm{d}\overline{q})\mathrm{d}p.$$
Similarly, if $\tilde{X}_{0}\sim\tilde{\mu}$, 
$$\tilde{X}_{\tilde{\tau}_\partial}\sim\frac{1}{\lambda_0+d\gamma}\left\vert p\cdot n(\overline{q})\right\vert\phi(\overline{q},p)\sigma_{\partial \mathcal{O}}(\mathrm{d}\overline{q})\mathrm{d}p.$$
\end{corollary}
% Furthermore, the exit event from $D$, starting from the respective QSD of the processes $(X_t)_{t\geq0}$ and $(\tilde{X}_t)_{t\geq0}$, satisfies the general properties stated in Proposition~\ref{exit event general case}. In addition, we show here that the law of the first exit point on $\partial D $ writes as follows. 

% \begin{proposition}[Law of the exit point] 
% Let the assumptions of Theorem~\ref{density qsd thm} hold. Assume that $X_0$ (resp. $\tilde{X}_0$) is distributed according to the QSD $\mu$ (resp. $\tilde{\mu}$) on $D$, then 
% \begin{enumerate}[label=(\roman*),ref=\roman*]
%     \item $X_{\tau_\partial}$ is distributed according to the density\quad$\frac{1}{\lambda_0}p\cdot n(q)\psi(q,p)\sigma_{\partial \mathcal{O}}(\mathrm{d}q)\mathrm{d}p$,
%     \item $\tilde{X}_{\tilde{\tau}_\partial}$ is distributed according to the density\quad$-\frac{1}{\lambda_0+d\gamma}p\cdot n(q)\phi(q,p)\sigma_{\partial \mathcal{O}}(\mathrm{d}q)\mathrm{d}p$,
% \end{enumerate} 
% where $\sigma_{\partial \mathcal{O}}(\mathrm{d}q)$ is the surface measure on $\partial \mathcal{O}$. 
% \end{proposition}

Furthermore, the densities of the QSD are the unique classical solutions to an eigenvalue problem.
\begin{theorem}[Spectral interpretation of the QSD]\label{vp droite thm} Under Assumptions~\ref{hyp F1 qsd} and~\ref{hyp O qsd}, there exists a unique couple $(\lambda,\eta)$ (resp. $(\lambda^*,\eta^*)$), up to a multiplicative constant on $\eta$ (resp. $\eta^*$) such that  $\eta\in\mathcal{C}^{2}(D)\cap\mathcal{C}^b(D\cup\Gamma^+)$ (resp. $\eta^*\in\mathcal{C}^{2}(D)\cap\mathcal{C}^b(D\cup\Gamma^-)$) is a non-zero, non-negative classical solution to the following problem
\begin{equation}\label{vp droite gauche}
  \left\{
    \begin{aligned}
        \mathcal{L}\eta(x) &=-\lambda \eta(x)&&\quad x\in D,\\
        \eta(x) &=0 &&\quad x\in \Gamma^+ ,
    \end{aligned}
\right.\qquad\text{resp.}\qquad
  \left\{\begin{aligned}
        \mathcal{L}^*\eta^*(x) &=-\lambda^* \eta^*(x)&&\quad x\in D,\\
        \eta^*(x) &=0 &&\quad x\in \Gamma^- .
    \end{aligned} \right. 
\end{equation}  
Moreover, $\eta \in\mathrm{L}^1(D)$, $\lambda=\lambda_0$ and $\frac{\eta}{\int_D \eta}=\phi$ (resp. $\eta^* \in\mathrm{L}^1(D)$, $\lambda^*=\lambda_0$ and $\frac{\eta^*}{\int_D \eta^*}=\psi$).
\end{theorem}
\begin{remark}
In particular, it follows from the expression of the spectral radii in Theorem~\ref{properties phi} that $\lambda_0$ is the smallest eigenvalue associated with the operators $-\mathcal{L}$ and $-\mathcal{L}^*$.
\end{remark}
Last, we are able to obtain the following long-time asymptotics of the operator $P^D_t$ on the Banach space $\mathrm{L}^\infty(D)$.

\begin{theorem}[Long-time asymptotics]\label{spectral decomp} Let Assumptions~\ref{hyp F1 qsd} and~\ref{hyp O qsd} hold. Let
$\alpha^*$ be defined by
\begin{equation}\label{eq:alpha star}
    \mathrm{e}^{-(\lambda_0+\alpha^*)}:=\sup_{z\in\sigma(P^D_1)\setminus\{\mathrm{e}^{-\lambda_0}\}} |z|.
\end{equation}
Then $\alpha^* \in (0,+\infty]$, and for all $\alpha\in[0,\alpha^*)$, there exists $C_\alpha>0$  such that for all $t\geq0$, for all $f\in\mathrm{L}^\infty(D)$,  
\begin{equation}\label{long time cv P^D_t}
    \left\Vert P^D_tf-\mathrm{e}^{-\lambda_0 t}\frac{\phi\otimes \psi(f)}{\int_D \phi \psi}\right\Vert_\infty\leq C_\alpha\Vert f\Vert_{\mathrm{L}^\infty(D)}\mathrm{e}^{-(\lambda_0+\alpha) t}, 
\end{equation} 
where the tensor product $\phi\otimes \psi$ is defined by: for any $f\in\mathrm{L}^\infty(D)$, $\phi\otimes \psi(f)=(\int_D \psi f)\phi$. 
\end{theorem}

\begin{remark}
In particular, for $f$ constant equal to $1$, Theorem~\ref{spectral decomp} ensures the existence of a constant $C>0$ such that for all $x\in D$, $t>0$,
$$\mathbb{P}_x(\tau_\partial>t)\leq C\mathrm{e}^{-\lambda_0 t}.$$ This ensures in particular that for any $\lambda\in[0,\lambda_0)$, 
$$\sup_{x\in D}\mathbb{E}_x\left[\mathrm{e}^{\lambda\tau_\partial}\right]<\infty.$$ 
\end{remark}

\begin{remark}
  A similar result can be deduced in the Banach space $L^p(D)$ for any $p\geq1$ using the compactness property of $P^D_t$ in $L^p(D)$ obtained in Theorem~\ref{compactness}. 
\end{remark}
We deduce from the previous result that the QSD $\mu$ attracts all probability measures $\theta$ on $D$ at an exponential rate.
\begin{theorem}[Convergence to the QSD in total variation]\label{cv semigroupe conditionnel} Under the assumptions of Theorem~\ref{spectral decomp}, for all $\alpha\in[0,\alpha^*)$, there exists $C'_\alpha>0$ such that, for all $t\geq0$, for any probability measure $\theta$ on $D$, $\mathbb{P}_{\theta}(\tau_\partial>t)>0$, and
\begin{equation}\label{semigroup conditionnel ineq}
    \left\|\mathbb{P}_\theta\left(X_t \in \cdot | \tau_\partial > t\right) - \mu\right\|_{TV} \leq \frac{C'_\alpha}{\int_D \phi \mathrm{d}\theta} \mathrm{e}^{-\alpha t},
\end{equation}
where $\|\cdot\|_{TV}$ denotes the total-variation norm on the space of bounded signed measures on $\mathbb{R}^{2d}$.
\end{theorem} 
\begin{remark}
The convergence speed and the prefactor are similar to what is obtained in~\cite[Chapter 4]{RamPHD} using the results by Champagnat and Villemonais in~\cite{V}.  
\end{remark}
\begin{remark}
Similar statements for Theorems~\ref{spectral decomp} and~\ref{cv semigroupe conditionnel} can be obtained for the adjoint Langevin process $(\Tilde{X}_t)_{t\geq0}$~\eqref{Langevin adjoint qsd} following the exact same proofs and using the equality in Theorem~\ref{duality thm qsd} satisfied by its transition density $\tilde{\mathrm{p}}^D_t(x,y)$. 
\end{remark}

The proofs of the theorems above are provided in the next Section.

\section{Proofs}\label{section results qsd} 

This section is organized  as follows. In Section~\ref{sec:compact}, we prove Theorem~\ref{compactness} and we prove Theorem~\ref{properties phi} in Section~\ref{section KreinRutman}. In Section~\ref{section qsd spectral interp} we obtain Theorems~\ref{density qsd thm},~\ref{prop: exit event qsd} and~\ref{vp droite thm} and finally, Section~\ref{section long time semigroupe} is devoted to the proofs of Theorems~\ref{spectral decomp} and~\ref{cv semigroupe conditionnel}.

\subsection{Proof of Theorem~\ref{compactness}}\label{sec:compact}
We prove in this section the results of Theorem~\ref{compactness} and in particular the compactness of the semigroup~$(P^D_t)_{t\geq0}$. The crucial ingredient is the Gaussian upper-bound satisfied by the transition density $\mathrm{p}_t^D$ (see Theorem~\ref{thm density intro qsd}): for any $\alpha\in (0,1)$ and for $t>0$, there exists $C_{t,t}>0$ such that for all $x,y\in D$,
\begin{equation}\label{Gaussian ub} 
    \mathrm{p}^D_t(x,y)\leq C_{t,t} \widehat{\mathrm{p}}^{(\alpha)}_{t}(x,y),
\end{equation}
where $\widehat{\mathrm{p}}^{(\alpha)}_{t}(x,y)$ is the transition density
of the Gaussian process
 $(\widehat{q}^{(\alpha)}_t, \widehat{p}^{(\alpha)}_t)_{t \geq 0}$ defined in~\eqref{eq:processus alpha}.  
 
 Let $\Phi_1,\Phi_2$ be the following positive continuous functions on $\mathbb{R}$:
\begin{equation}\label{Phi_1}
\Phi_1:\rho\in\mathbb{R}\mapsto\begin{cases}
  \frac{1-\mathrm{e}^{-\rho}}{\rho}&\text{if $\rho\neq0$,}\\ 
  1 & \text{if $\rho=0$,} 
\end{cases}
\end{equation}
\begin{equation}\label{Phi_2}
\Phi_2:\rho\in\mathbb{R}\mapsto\begin{cases}
  \frac{3}{2 \rho^3}\left[2\rho-3+4 \mathrm{e}^{-\rho}-\mathrm{e}^{-2\rho}\right] & \text{if $\rho\neq0$,}\\
  1 &\text{if $\rho=0$.}
\end{cases}
\end{equation}
 
 One can show, see~\cite[Section 5.1]{LelRamRey}, that for all $t \geq 0$ and $\alpha\in(0,1]$, the vector $(\widehat{q}_t^{(\alpha)},\widehat{p}_t^{(\alpha)})$ admits the following law under $\mathbb{P}_{(q,p)}$
\begin{equation}\label{loi gaussienne}
\begin{pmatrix}
\widehat{q}^{(\alpha)}_t  \\
\widehat{p}^{(\alpha)}_t
\end{pmatrix}
  \sim  \mathcal{N}_{2d}\left(
\begin{matrix} 
\begin{pmatrix}
m_q(t)  \\
m_p(t) 
\end{pmatrix},
\displaystyle{\frac{C(t)}{\alpha}}
\end{matrix} 
\right)  , 
\end{equation}
where the mean vector is $$m_q(t):=q+t p \Phi_1(\gamma t),\qquad m_p(t):=p \mathrm{e}^{-\gamma  t},$$
and the matrix $C(t)$ is
$$C(t):= \begin{pmatrix}
c_{qq}(t)  I_d & c_{qp}(t)  I_d \\
c_{qp}(t)  I_d & c_{pp}(t)  I_d 
\end{pmatrix}, 
$$
where $I_d$ is the identity matrix in $\mathbb{R}^{d\times d}$ and
\begin{equation}\label{coeff cov}
    c_{qq}(t):=\frac{\sigma^2 t^3}{3} \Phi_2(\gamma t),\qquad c_{qp}(t):=\frac{\sigma^2 t^2}{2} \Phi_1(\gamma t)^2,\qquad c_{pp}(t):=\sigma^2t \Phi_1(2\gamma t)  . 
\end{equation} 
The determinant of the covariance matrix $\frac{C(t)}{\alpha}$ is $\mathrm{det}\left(\frac{C(t)}{\alpha}\right)=\left(\frac{\sigma^4 t^4}{12\alpha^2} \phi(\gamma t)\right)^d$ where $\phi$ is the positive continuous function defined by 
\begin{equation}\label{expr phi}
\phi:\rho\in\mathbb{R}\mapsto 4 \Phi_2(\rho) \Phi_1(2 \rho)-3 \Phi_1(\rho)^4=\begin{cases}
    \frac{6 (1-\mathrm{e}^{-\rho})}{\rho^4}\left[-2+\rho+(2+\rho)\mathrm{e}^{-\rho}\right]&\text{if $\rho\neq0$,}\\
   1 &\text{if $\rho=0$.}
\end{cases}
\end{equation}

Using the Gaussian upper-bound~\eqref{Gaussian ub} let us now prove that $\mathrm{p}^D_t\in\mathrm{L}^\infty(D\times D)\cap\mathrm{L}^1(D\times D)$.
\begin{lemma}[$\mathrm{p}^D_t\in\mathrm{L}^\infty(D\times D)\cap\mathrm{L}^1(D\times D)$]\label{hilbert schmidt lemma} Let Assumptions \ref{hyp O qsd} and \ref{hyp F1 qsd} hold.
For all $t>0$, $\mathrm{p}_t^D$ is bounded on $D\times D$ and
\begin{equation}\label{hilbert schmidt}
    \iint_{D\times D} \mathrm{p}_t^D(x,y) \mathrm{d}x \mathrm{d}y<\infty  . 
\end{equation} 
\end{lemma} 
\begin{proof}
Let $t>0$ and $\alpha\in(0,1)$. By~\eqref{Gaussian ub}, for all $x,y\in D$,
$$\mathrm{p}^D_t(x,y)\leq\frac{C_{t,t}}{ \sqrt{(2 \pi)^{2d} \mathrm{det}\left(\frac{C(t)}{\alpha}\right)}}=\frac{C_{t,t}}{ \sqrt{(2 \pi)^{2d} \left(\frac{\sigma^4 t^4}{12\alpha^2} \phi(\gamma t)\right)^d}},$$
which ensures that $\mathrm{p}^D_t\in\mathrm{L}^\infty(D\times D)$. To obtain that $\mathrm{p}^D_t\in\mathrm{L}^1(D\times D)$, we prove that 
$$\iint_{D\times D} \widehat{\mathrm{p}}^{(\alpha)}_t(x,y) \mathrm{d}x\mathrm{d}y<\infty.$$ 

Let us first integrate $\widehat{\mathrm{p}}^{(\alpha)}_t((q,p),(q',p'))$ with respect to $p,p'\in\mathbb{R}^d$ using
Fubini-Tonelli's theorem. Since
$\widehat{\mathrm{p}}^{(\alpha)}_t((q,p),(q',p'))$ is the transition
density of the Gaussian process $(\widehat{q}^{(\alpha)}_t, \widehat{p}^{(\alpha)}_t)_{t \geq 0}$, one can obtain an explicit
expression of
$\int_{\mathbb{R}^d}\widehat{\mathrm{p}}^{(\alpha)}_t((q,p),(q',p'))
\mathrm{d}p'$ which corresponds to the marginal density of $\widehat{q}^{(\alpha)}_t$ under $\mathbb{P}_{(q,p)}$. Since, under $\mathbb{P}_{(q,p)}$, 
\begin{equation*}
  \widehat{q}^{(\alpha)}_t \sim \mathcal{N}_d\left(q + tp\Phi_1(\gamma t), \frac{c_{qq}(t)}{\alpha}I_d\right),\quad \frac{c_{qq}(t)}{\alpha} = \frac{\sigma^2 t^3}{3\alpha}\Phi_2(\gamma t),
\end{equation*}
so that  
 \begin{align*}
   \int_{\mathbb{R}^d}\widehat{\mathrm{p}}^{(\alpha)}_t((q,p),(q',p')) \mathrm{d}p'&=\frac{(3\alpha)^{d/2}}{\left(2 \pi\sigma^2 t^3 \Phi_2(\gamma t)\right)^{d/2}} \mathrm{e}^{-\frac{3 \alpha}{2 \sigma^2 t^3 \Phi_2(\gamma t)}\left\vert q'-q-t p \Phi_1(\gamma t) \right\vert^2} 
\end{align*}
where $\Phi_1$ and $\Phi_2$ are defined in \eqref{Phi_1} and \eqref{Phi_2}. Then,  
\begin{align*}
   \int_{\mathbb{R}^d}\int_{\mathbb{R}^d}\widehat{\mathrm{p}}^{(\alpha)}_t((q,p),(q',p')) \mathrm{d}p \mathrm{d}p'&=\frac{(3\alpha)^{d/2}}{\left(2 \pi\sigma^2 t^3 \Phi_2(\gamma t)\right)^{d/2}}\int_{\mathbb{R}^d}\mathrm{e}^{-\frac{3 \alpha}{2 \sigma^2 t^3 \Phi_2(\gamma t)}\left\vert q'-q-t p \Phi_1(\gamma t) \right\vert^2} \mathrm{d}p\\
   &=\frac{1}{t^d\Phi_1(\gamma t)^d}.
\end{align*}
Consequently,
$$\iint_{D\times D} \widehat{\mathrm{p}}^{(\alpha)}_t(x,y) \mathrm{d}x\mathrm{d}y=\frac{\vert\mathcal{O}\vert^2}{t^d\Phi_1(\gamma t)^d},$$ which concludes the proof.
\end{proof}
\begin{remark}\label{p^D_t in any L^q}
This ensures in particular that $\mathrm{p}^D_t\in\mathrm{L}^q(D\times D)$ for any $q\in[1,+\infty]$ and $t>0$.
\end{remark}
Let us now prove Theorem \ref{compactness}. 
\begin{proof}[Proof of Theorem \ref{compactness}]

It follows from Remark~\ref{p^D_t in any L^q} that for all $t>0$, $$\iint_{D\times D} \mathrm{p}_t^D(x,y)^2 \mathrm{d}x \mathrm{d}y<\infty .$$ Therefore, by Theorems VI.22 and VI.23 in \cite{RS}, the operator $P^D_t$ is a compact operator from $\mathrm{L}^2(D)$ to $\mathrm{L}^2(D)$. 

Let $s>0$. In \textbf{Step~1}, we show that $P_s^D$ maps $\mathrm{L}^p(D)$, $p \in [1,+\infty)$, continuously into $\mathrm{L}^\infty(D)$. In \textbf{Step~2}, we show that $P_s^D$ maps $\mathrm{L}^\infty(D)$ continuously into $\mathrm{L}^q(D)$, $q \in [1,+\infty)$. In \textbf{Step~3}, we show that $P_s^D$ maps $\mathrm{L}^\infty(D)$ continuously into $\mathcal{C}^b(\overline{D})$, and that $P_s^D$ satisfies a semigroup property. We conclude the proof in \textbf{Step~4}.

\medskip \noindent \textbf{Step 1}. Let $s>0$, $p\geq1$. Let us prove that $P^D_s$ maps $\mathrm{L}^p(D)$ continuously into $\mathrm{L}^\infty(D)$. Recall that by Lemma \ref{hilbert schmidt lemma}, $\|\mathrm{p}^D_s\|_{\mathrm{L}^\infty(D \times D)} < +\infty$. Therefore if $p=1$, then for any $f \in \mathrm{L}^1(D)$ we have
\begin{equation*}
  \|P^D_s f\|_{\mathrm{L}^\infty(D)} \leq \|\mathrm{p}^D_s\|_{\mathrm{L}^\infty(D \times D)} \|f\|_{\mathrm{L}^1(D)},
\end{equation*}
while if $p \in (1,+\infty)$, then for any $f \in \mathrm{L}^p(D)$ and $x \in D$, letting $q \in (1,+\infty)$ be such that $1/p+1/q=1$, we get by H\"older's inequality
\begin{align*}
  \left|P^D_sf(x)\right| &\leq \|\mathrm{p}^D_s(x,\cdot)\|_{\mathrm{L}^q(D)} \|f\|_{\mathrm{L}^p(D)}\\
  &\leq \left(\|\mathrm{p}^D_s\|_{\mathrm{L}^\infty(D \times D)}^{q-1}\mathbb{P}(\tau^x_\partial>s)\right)^{1/q} \|f\|_{\mathrm{L}^p(D)}\\
  &\leq \|\mathrm{p}^D_s\|_{\mathrm{L}^\infty(D \times D)}^{\frac{q-1}{q}} \|f\|_{\mathrm{L}^p(D)},
\end{align*}
which yields $\|P^D_s f\|_{\mathrm{L}^\infty(D)} \leq \|\mathrm{p}^D_s\|_{\mathrm{L}^\infty(D \times D)}^{\frac{q-1}{q}}  \|f\|_{\mathrm{L}^p(D)}$.

\medskip \noindent \textbf{Step 2}. Let $s>0$, $q\geq1$. Let us prove that $P^D_s$ maps $\mathrm{L}^\infty(D)$ continuously into $\mathrm{L}^q(D)$. Let $f\in\mathrm{L}^\infty(D)$. For $x\in D$, one has that
\begin{align*}
   \bigg\vert P^D_s f(x)\bigg\vert^q&=\bigg\vert\int_D \mathrm{p}_s^D(x,y) f(y) \mathrm{d}y\bigg\vert^q\\
   &\leq \mathbb{P}(\tau^x_\partial>s)^{q}\Vert f\Vert_{\mathrm{L}^\infty(D)}^q\\
   &\leq \mathbb{P}(\tau^x_\partial>s)\Vert f\Vert_{\mathrm{L}^\infty(D)}^q.
\end{align*}
Therefore, using Lemma \ref{hilbert schmidt lemma} we get $\Vert P^D_sf\Vert_{\mathrm{L}^q(D)}\leq \|\mathrm{p}^D_s\|_{\mathrm{L}^1(D \times D)}^{1/q}\Vert f\Vert_{\mathrm{L}^\infty(D)}$.

\medskip \noindent \textbf{Step 3}. We first deduce from the Markov property that
\begin{equation*}
  P^D_{t+s}f(x) = \mathbb{E}_x\left[f(X_{t+s})\mathbb{1}_{\tau_\partial>t+s}\right] = \mathbb{E}_x\left[\mathbb{1}_{\tau_\partial>t}P^D_sf(X_t)\right] = P^D_t P^D_s f,
\end{equation*}
which together with the obvious observation that $P^D_0f=f$, shows that $(P^D_t)_{t \geq 0}$ is a semigroup on $\mathrm{L}^p(D)$ for any $p \in [1,+\infty]$ by \textbf{Steps 1} and \textbf{2}. Let $s>0$ and $f \in \mathrm{L}^\infty(D)$. Since for all $x\in D$,
\begin{align*}
    P^D_s f(x)=\int_D \mathrm{p}_{s/2}^D(x,y) P^D_{s/2}f(y) \mathrm{d}y,
\end{align*}
and $P^D_{s/2}$ maps $\mathrm{L}^\infty(D)$ into $\mathrm{L}^1(D)$ by \textbf{Step 2}, then it is an immediate application of the dominated convergence theorem along with the continuity and boundedness of $x\in\overline{D}\mapsto\mathrm{p}_{s/2}^D(x,y)$ in Theorem~\ref{thm density intro qsd} that  $P^D_sf \in \mathcal{C}^b(\overline{D})$. Besides, one has obviously 
\begin{equation*}
  \|P^D_sf\|_{\mathcal{C}^b(\overline{D})} := \|P^D_sf\|_\infty \leq \|f\|_{L^\infty(D)},
\end{equation*}
so that $P^D_s$ maps $\mathrm{L}^\infty(D)$ into $\mathcal{C}^b(\overline{D})$ continuously. Notice that  $(P^D_t)_{t \geq 0}$ is thus also a semigroup on $\mathcal{C}^b(\overline{D})$.

\medskip \noindent \textbf{Step 4}. In order to study compactness, we shall use repeatedly the fact that the composition of a continuous operator with a compact operator is a compact operator.

Let $p \in [1,+\infty]$ and $t>0$. Writing $P^D_t = P^D_{t/3} P^D_{t/3} P^D_{t/3}$, using the continuity of the mappings $P^D_{t/3} : \mathrm{L}^p(D) \to \mathrm{L}^2(D)$ and $P^D_{t/3} : \mathrm{L}^2(D) \to \mathrm{L}^p(D)$, and the compactness of $P^D_{t/3} : \mathrm{L}^2(D) \to \mathrm{L}^2(D)$, we obtain that $P^D_t$ is a compact operator from $\mathrm{L}^p(D)$ to $\mathrm{L}^p(D)$.

Similarly, writing $P^D_t = \iota P^D_{t/2}P^D_{t/2}$, where $\iota$ is the injection from $\mathcal{C}^b(\overline{D})$ to $\mathrm{L}^\infty(D)$, using the continuity of the operators $\iota$ and $P_{t/2} : \mathrm{L}^\infty(D) \to \mathcal{C}^b(\overline{D})$, as well as the compactness of the operator $P^D_{t/2} : \mathrm{L}^\infty(D) \to \mathrm{L}^\infty(D)$ that we have just proven, we conclude that $P^D_t$ is a compact operator from $\mathcal{C}^b(\overline{D})$ to $\mathcal{C}^b(\overline{D})$.
\end{proof}  
\subsection{Proof of Theorem~\ref{properties phi}}
\label{section KreinRutman} 
Using the compactness properties obtained in Theorem~\ref{compactness} and Remark~\ref{semigroupe adjoint},  we apply the Krein-Rutman theorem to the operators $P^D_t$ and $\Tilde{P}^D_t$ defined in~\eqref{def semigroupe} and~\eqref{def semigroupe adjoint}. In order to do that let us first recall an important property satisfied by the spectral radius of a bounded operator, see~\cite[p. 192, Theorem VI.6]{RS}. 
\begin{proposition}[Gelfand's formula]\label{spectral radius} 
Let $T$ be a bounded real operator on $\mathcal{C}^b(\overline{D})$. One has that
$$r(T)=\lim_{n\rightarrow\infty}\vertiii{T^n}_{\mathcal{C}^b(\overline{D})}^{1/n}.$$ 
\end{proposition}
 The Krein-Rutman theorem, recalled below, basically states that under some conditions on the bounded operator $T$, the spectral radius $r(T)$ is also an eigenvalue of the operator $T$. The following version of the Krein-Rutman theorem can be found in~\cite[p. 313]{SW}.
 \begin{theorem}[Krein-Rutman] \label{KreinRutman} Let $K\subset \mathcal{C}^b(\overline{D})$ be a convex cone such that the set $\{f-g:f,g\in K\}$ is dense in $\mathcal{C}^b(\overline{D})$. Let $T:\mathcal{C}^b(\overline{D})\mapsto \mathcal{C}^b(\overline{D})$ be a compact operator such that $T(K)\subset K$, and assume that its spectral radius $r(T)$ is strictly positive. Then, $r(T)$ is an eigenvalue of $T$ with an eigenvector $u$ in $K\setminus\{0\}$ such that $T(u)=r(T)u$. 
\end{theorem}
Let $K$ be the following convex cone of $\mathcal{C}^b(\overline{D})$, \begin{equation}\label{cone K}
    K:=\{f\in\mathcal{C}^b(\overline{D}) : f\geq0\}.
\end{equation} The density of $\{f-g:f,g\in K\}$ in $\mathcal{C}^b(\overline{D})$ is immediate. Our goal now is to apply the Krein-Rutman theorem above to the compact operators  $P^D_t$ and $\Tilde{P}^D_t$, for $t>0$, on the cone $K$. In order to do that we need to prove the positivity of the spectral radii $r(P^D_t)$ and $r(\Tilde{P}^D_t)$ for $t>0$.
\begin{proposition}[Spectral radius positivity]\label{prop spectral radius} Under Assumptions~\ref{hyp F1 qsd} and~\ref{hyp O qsd}, for all $t>0$, $r(P^D_t)>0$ and $r(\Tilde{P}^D_t)>0$.  
\end{proposition}
\begin{proof}
Let $t>0$. We prove here that $r(P^D_t)>0$ which relies merely on the positivity of its transition density $\mathrm{p}^D_t(\cdot,\cdot)$ on $D\times D$, see Theorem~\ref{thm density intro qsd}. Besides, the positivity of the transition density $\tilde{\mathrm{p}}^D_t(\cdot,\cdot)$ easily follows from the equality in Theorem~\ref{duality thm qsd}, therefore the exact same proof applies to $\Tilde{P}^D_t$ and ensures that $r(\Tilde{P}^D_t)>0$. As a result, we omit here the case of $\Tilde{P}^D_t$ to avoid repetition. 

Following Proposition~\ref{spectral radius}, it is sufficient to prove that there exists a constant $\beta>0$ such that for all $n\geq1$, $\vertiii{(P^D_t)^n}_{\mathcal{C}^b(\overline{D})}^{1/n}\geq \beta$. Let $C\subset D$ be a compact set with positive Lebesgue measure, i.e. $\vert C\vert>0$. It follows from Theorem~\ref{thm density intro qsd} that $P^D_t$ admits a smooth transition density $\mathrm{p}^D_t$, which is positive on $C\times C$. Therefore, there exists $\alpha>0$ such that for all $x,y\in C$, 
$\mathrm{p}^D_t(x,y)\geq\alpha$. Besides, for $n\geq1$, $$\vertiii{\left(P^D_t\right)^n}_{\mathcal{C}^b(\overline{D})}\geq \left\Vert\left(P^D_t\right)^n\mathbb{1}_D\right\Vert_\infty .$$ Moreover, for all $x\in C$,
\begin{align*}
    \left\Vert\left(P^D_t\right)^n\mathbb{1}_D\right\Vert_\infty&\geq \left(P^D_t\right)^n\mathbb{1}_D(x)\\
    &= \int_{D^n}\mathrm{p}^D_t(x,y_1)\dots\mathrm{p}^D_t(y_{n-1},y_n) \mathrm{d}y_1\dots \mathrm{d}y_n\\
    &\geq \int_{C^n}\mathrm{p}^D_t(x,y_1)\dots\mathrm{p}^D_t(y_{n-1},y_n) \mathrm{d}y_1\dots \mathrm{d}y_n\\
    &\geq (\alpha \vert C\vert)^n .
\end{align*}
Consequently, for all $n\geq1$, $$\vertiii{\left(P^D_t\right)^n}_{\mathcal{C}^b(\overline{D})}^{1/n}\geq\alpha \vert C\vert>0 ,$$ which concludes the proof.
\end{proof} 

We now apply Theorem~\ref{KreinRutman} to the operators $P^D_t$ and $ \Tilde{P}^D_t$ on the cone $K$ to obtain that $r(P^D_t)$ and $r(\Tilde{P}^D_t)$ are eigenvalues of their respective operators. In addition, these eigenvalues are shown to be simple. The following proof is inspired from unpublished lecture notes by P. Collet. 
\begin{proposition}[Consequence of the Krein-Rutman  theorem]\label{spectral radius lemma} Let Assumptions~\ref{hyp F1 qsd} and~\ref{hyp O qsd} hold. For all $t>0$, $r(P^D_t)$ (resp. $r(\Tilde{P}^D_t)$) is a simple eigenvalue of the operator $P^D_t$ (resp. $\Tilde{P}^D_t$) with eigenspace generated by an element $\phi_t$ (resp. $\psi_t$) of $K\cap\mathrm{L}^1(D)$ such that $\phi_t>0$ (resp. $\psi_t>0$) on $D$. 
Furthermore, $r(\Tilde{P}^D_t)=r(P^D_t) \mathrm{e}^{-d\gamma t}$.  
\end{proposition}
\begin{proof}
Let $t>0$. The compactness of $P^D_t$ and $\Tilde{P}^D_t$ on $\mathcal{C}^b(\overline{D})$ follows from Theorem~\ref{compactness} and Remark~\ref{semigroupe adjoint}. Besides, the cone $K$ defined in~\eqref{cone K} evidently satisfies the assumptions of Theorem~\ref{KreinRutman}. Therefore, by Theorem~\ref{KreinRutman} and Proposition~\ref{prop spectral radius} we obtain the existence of $\phi_t,\psi_t\in K\setminus\{0\}$ such that
\begin{equation}\label{eq phi psi}
    P^D_t\phi_t=r(P^D_t)\phi_t,\qquad  \Tilde{P}^D_t\psi_t=r(\Tilde{P}^D_t)\psi_t.
\end{equation}

\medskip \noindent \textbf{Step 1}. Let us prove that 
\begin{equation}\label{rayon spectral}
    r(\Tilde{P}^D_t)=r(P^D_t) \mathrm{e}^{-d\gamma t} 
\end{equation}
by computing the integral $\int_D\psi_t(x) P^D_t\phi_t(x) \mathrm{d}x$ in two different ways. On the one hand, by~\eqref{eq phi psi},
\begin{equation}\label{spectral radius eq 1}
    \int_D\psi_t(x) P^D_t\phi_t(x) \mathrm{d}x=r(P^D_t) \int_D\psi_t(x) \phi_t(x) \mathrm{d}x .
\end{equation}
On the other hand, using Fubini-Tonelli's theorem, Theorem~\ref{duality thm qsd} and~\eqref{eq phi psi} again,
\begin{align*}
\int_D\psi_t(x) P^D_t\phi_t(x) \mathrm{d}x&=\mathrm{e}^{d\gamma t}\int_D\phi_t(x) \Tilde{P}^D_t\psi_t(x) \mathrm{d}x \\
&=\mathrm{e}^{d\gamma t}r(\Tilde{P}^D_t)\int_D\phi_t(x)\psi_t(x) \mathrm{d}x.\numberthis\label{spectral radius eq 2}
\end{align*} 
Let us now prove that $\int_D\phi_t(x)\psi_t(x) \mathrm{d}x\in(0,\infty)$. First, for $x\in D$, $r(P^D_t) \phi_t(x)=P^D_t\phi_t(x)=\int_D\mathrm{p}^D_t(x,y) \phi_t(y) \mathrm{d}y>0$ since $\phi_t\in K\setminus\{0\}$ and $\mathrm{p}^D_t>0$ on $D\times D$ by Theorem~\ref{thm density intro qsd}. Therefore, $\phi_t>0$ on $D$. Likewise $\psi_t>0$ on $D$\ so that $\int_D\phi_t(x)\psi_t(x) \mathrm{d}x>0$. Second, using the boundedness of $\phi_t$ along with the double integral estimate from Lemma~\ref{hilbert schmidt lemma} in the left equality in~\eqref{eq phi psi}, one obtains that $\phi_t\in \mathrm{L}^1(D)$. Using Theorem~\ref{duality thm qsd} one also has that $\psi_t\in \mathrm{L}^1(D)$. In particular, since $\phi_t$ and $\psi_t$ are in $\mathrm{L}^\infty(D)$, this yields that $\int_D\phi_t(x)\psi_t(x)\mathrm{d}x<\infty$. As a result, the equalities~\eqref{spectral radius eq 1} and~\eqref{spectral radius eq 2} yield~\eqref{rayon spectral}. 

\medskip \noindent \textbf{Step 2}. Let us prove that every real-valued eigenvector of $P^D_t$ associated with the eigenvalue $r(P^D_t)$ has a constant sign. Assume that there exists $h_t\in\mathcal{C}^b(\overline{D})$ such that $P^D_th_t=r(P^D_t)h_t$ and $h_t$ changes sign on $D$. Then, by the positivity of $\mathrm{p}^D_t$ one has for $x\in D$,
\begin{align*}
    r(P^D_t)\vert h_t(x)\vert&=\vert P^D_th_t(x)\vert\\
    &=\Big\vert \int_D\mathrm{p}^D_t(x,y)h_t(y)\mathrm{d}y\Big\vert\\
    &< \int_D\mathrm{p}^D_t(x,y)\vert h_t(y)\vert\mathrm{d}y=P^D_t\vert h_t\vert(x) .
\end{align*}
As a result, since $\psi_t>0$ on $D$, by Theorem~\ref{duality thm qsd} one has that
\begin{align*}
    r(P^D_t)\int_D\psi_t(x) \vert h_t(x)\vert\mathrm{d}x&<\int_D\psi_t(x) P^D_t\vert h_t\vert(x)\mathrm{d}x\\
    &=\mathrm{e}^{d\gamma t}\int_D\Tilde{P}^D_t\psi_t(x) \vert h_t(x)\vert\mathrm{d}x\\
    &=r(P^D_t)\int_D\psi_t(x) \vert h_t(x)\vert\mathrm{d}x,
\end{align*}
by~\eqref{rayon spectral}, which leads to a contradiction. Therefore, $h_t$ has a constant sign on $D$. 

\medskip \noindent \textbf{Step 3}. Let $h_t$ be a real-valued eigenvector of $P^D_t$ associated with the eigenvalue $r(P^D_t)$, let us prove that $h_t\in\mathrm{Span}(\phi_t)$. Up to changing $h_t$ to $-h_t$, we can assume that $h_t\in K$ by \textbf{Step 2}. Let us define for $x\in D$, $$\Tilde{h}_t(x):=\frac{h_t(x)}{\int_D \psi_t(y)h_t(y)\mathrm{d}y},\qquad\Tilde{\phi}_t(x):=\frac{\phi_t(x)}{\int_D \psi_t(y)\phi_t(y)\mathrm{d}y},$$ so that
\begin{equation}\label{vp tilde}
    \int_D\Tilde{h}_t(x)\psi_t(x)\mathrm{d}x=\int_D\Tilde{\phi}_t(x)\psi_t(x)\mathrm{d}x=1.
\end{equation} 
Notice that, $\Tilde{\phi}_t-\Tilde{h}_t$ is an eigenvector of $P^D_t$ with eigenvalue $r(P^D_t)$, therefore it has a constant sign. By~\eqref{vp tilde}, one concludes that necessarily  $\Tilde{\phi}_t-\Tilde{h}_t=0$ on $D$ since $\psi_t>0$ on $D$. Hence $h_t\in\mathrm{Span}(\phi_t)$ and $r(P^D_t)$ is a simple eigenvalue. 

\medskip \noindent \textbf{Step 4}. Applying this time \textbf{Step 2} and \textbf{Step 3} to the operator $\Tilde{P}^D_t$, one also obtains that $r(\Tilde{P}^D_t)$ is also a simple eigenvalue. This concludes the proof of Proposition~\ref{spectral radius lemma}.
\end{proof}
To prove Theorem~\ref{properties phi}, we finally need the following technical lemma.
\begin{lemma}[High velocity exit event]\label{conv vitesse infinie} Under Assumptions~\ref{hyp F1 qsd} and~\ref{hyp O qsd},
$$\forall t>0,\qquad\sup_{q\in\mathcal{O}}\mathbb{P}_{(q,p)}(\tau_\partial>t)\underset{\vert p\vert\rightarrow\infty}{\longrightarrow}0.$$ 
\end{lemma}
\begin{proof}
Let $t>0$ and $\alpha\in(0,1)$. The Gaussian upper-bound~\eqref{Gaussian ub} ensures the existence of $C_{t,t}>0$ such that for all $x,y\in D$, $\mathrm{p}_t^D(x,y)\leq C_{t,t} \widehat{\mathrm{p}}^{(\alpha)}_t(x,y)$, where $\widehat{\mathrm{p}}^{(\alpha)}_t(x,y)$ is the transition density of the process $(\widehat{q}^{(\alpha)}_t, \widehat{p}^{(\alpha)}_t)_{t \geq 0}$ defined in~\eqref{eq:processus alpha}.  

Furthermore, for $x=(q,p)\in D$, the law of  $\widehat{q}^{(\alpha)}_t$ is Gaussian with mean $q+t p \Phi_1(\gamma t)$ and covariance matrix $\frac{\sigma^2 t^3}{3\alpha} \Phi_2(\gamma t)I_d$, where $\Phi_1$ are $\Phi_2$ are defined in~\eqref{Phi_1} and~\eqref{Phi_2}. Therefore, 
 \begin{align*}
   \int_{\mathbb{R}^d}\widehat{\mathrm{p}}^{(\alpha)}_t((q,p),(q',p')) \mathrm{d}p'&=\frac{(3\alpha)^{d/2}}{\left(2 \pi\sigma^2 t^3 \Phi_2(\gamma t)\right)^{d/2}} \mathrm{e}^{-\frac{3 \alpha t^2 \Phi_1(\gamma t)^2}{2 \sigma^2 t^3 \Phi_2(\gamma t)}\left\vert p-\frac{q'-q}{t \Phi_1(\gamma t)} \right\vert^2}. 
\end{align*}  

Let $\delta:=\sup_{q,q'\in\mathcal{O}}\vert q-q'\vert$ (which is finite since $\mathcal{O}$ is bounded), then for $t>0$ and $q,q'\in\mathcal{O}$,  
if $\vert p\vert\geq\frac{2 \delta}{t \Phi_1(\gamma t)}$ ($\Phi_1$ is positive),
$$\int_{\mathbb{R}^d}\widehat{\mathrm{p}}^{(\alpha)}_t((q,p),(q',p')) \mathrm{d}p'\leq\frac{(3\alpha)^{d/2}}{\left(2 \pi\sigma^2  t^3 \Phi_2(\gamma t)\right)^{d/2}} \mathrm{e}^{-\frac{3 \alpha t^2 \Phi_1(\gamma t)^2}{8\sigma^2 t^3 \Phi_2(\gamma t)}\left\vert p\right\vert^2}. $$
As a consequence,  
\begin{align*}
   \sup_{q \in \mathcal{O}}\mathbb{P}_{(q,p)}(\tau_\partial>t)&=\sup_{q \in \mathcal{O}}\int_{\mathcal{O}}\int_{\mathbb{R}^d}\mathrm{p}^D_t((q,p),(q',p')) \mathrm{d}p' \mathrm{d}q'\\
   &\leq \sup_{q \in \mathcal{O}} C_{t,t}\int_{\mathcal{O}}\int_{\mathbb{R}^d}\widehat{\mathrm{p}}^{(\alpha)}_t((q,p),(q',p')) \mathrm{d}p' \mathrm{d}q'\\
   &\leq C_{t,t}\frac{(3\alpha)^{d/2}\vert\mathcal{O}\vert}{\left(2 \pi\sigma^2 t^3 \Phi_2(\gamma t)\right)^{d/2}} \mathrm{e}^{-\frac{3 \alpha t^2 \Phi_1(\gamma t)^2}{8\sigma^2 t^3 \Phi_2(\gamma t)}\left\vert p\right\vert^2}\underset{\vert p\vert\rightarrow\infty}{\longrightarrow}0, 
\end{align*}
which concludes the proof.
\end{proof}
We are now in position to prove Theorem~\ref{properties phi}.
\begin{proof}[Proof of Theorem~\ref{properties phi}]
 For $t>0$, let $\phi_t$ (resp. $\psi_t$) be an eigenvector of $P^D_t$ (resp. $\Tilde{P}^D_t$) in $K\setminus\{0\}$ associated with the eigenvalue $r(P^D_t)$ (resp. $r(\Tilde{P}^D_t)$) and such that $\phi_t>0$ (resp. $\psi_t>0$) on $D$, whose existence is ensured by Proposition~\ref{spectral radius lemma}. We will prove Theorem~\ref{properties phi} for $\phi_t$ and $P^D_t$, but the exact same reasoning with the operator $\Tilde{P}^D_t$ instead yields the proof for $\psi_t$ and $\tilde{P}^D_t$. 
 
\medskip \noindent
 \textbf{Step 1}. Let us start by proving that $\mathrm{Span}(\phi_t , t>0)$ is a one-dimensional space generated by a function $\phi\in K$. This is the case if one can prove that for all $s,t>0$, $\phi_s\in\mathrm{Span}(\phi_t)$. 
 
 For $s,t>0$, $P^D_s\phi_s=r(P^D_s)\phi_s$. Furthermore, by the semigroup property satisfied by $(P^D_r)_{r\geq0}$, $$P^D_sP^D_t\phi_s=P^D_tP^D_s\phi_s=r(P^D_s)P^D_t\phi_s.$$ Since $r(P^D_s)$ is a simple eigenvalue of $P^D_s$ by Proposition~\ref{spectral radius lemma} then $P^D_t\phi_s\in\mathrm{Span}(\phi_s)$, i.e. there exists $\alpha_{t,s}>0$ such that $P^D_t\phi_s=\alpha_{t,s}\phi_s$.  
 
Let us prove that $\alpha_{t,s}=r(P^D_t)$. Consider the integral $\int_DP^D_t\phi_s(x)\psi_t(x)\mathrm{d}x$. One has that
 $$\int_DP^D_t\phi_s(x)\psi_t(x)\mathrm{d}x=\alpha_{t,s}\int_D\phi_s(x)\psi_t(x)\mathrm{d}x.$$
 Furthermore, Theorem~\ref{duality thm qsd} and Proposition~\ref{spectral radius lemma} also ensure that
 $$\int_DP^D_t\phi_s(x)\psi_t(x)\mathrm{d}x=\mathrm{e}^{d\gamma t}\int_D\phi_s(x)\Tilde{P}^D_t\psi_t(x)\mathrm{d}x=r(P^D_t)\int_D\phi_s(x)\psi_t(x)\mathrm{d}x.$$
Since $\phi_s,\psi_t$ are positive on $D$ and belong to $\mathcal{C}^b(\overline{D})\cap\mathrm{L}^1(D)$ then $\int_D\phi_s(x)\psi_t(x)\mathrm{d}x\in(0,\infty)$. Therefore, the equalities above ensure that $\alpha_{t,s}=r(P^D_t)$. In particular, this yields that $\phi_s\in\mathrm{Span}(\phi_t)$ since $r(P^D_t)$ is a simple eigenvalue for $P^D_t$. Let us now denote by $\phi\in K$ a function generating $\mathrm{Span}(\phi_t , t>0)$. 

\medskip \noindent \textbf{Step 2}. Let us now show that there exists $\lambda_{0}\geq0$ such that for all $t>0$, $r(P^D_t)=\mathrm{e}^{-\lambda_0t}$. For $s,t>0$, $P^D_{t+s}\phi=r(P^D_{t+s})\phi$. Besides,
\begin{align*}
    P^D_{t+s}\phi&=P^D_tP^D_s\phi\\
    &=r(P^D_t)r(P^D_s)\phi .
\end{align*}
Therefore, $r(P^D_{t+s})=r(P^D_t)r(P^D_s)$ since $\phi>0$. Since $\vertiii{P^D_s}_{\mathcal{C}^b(\overline{D})}\leq1$ for all $s>0$, then $r(P^D_s)\leq1$. As a result, for all $s,t>0$, $r(P^D_{t+s})\leq r(P^D_t)$. Consequently, the function $v:t>0\mapsto \log(r(P^D_t))$ is a decreasing function which satisfies the Cauchy equation:  
$$\forall s,t>0,\qquad v_{t+s}=v_t+v_s .$$
Classical results for Cauchy equations ensure that $v_t$ is linear. This implies that there exists $\lambda_0\geq0$ such that for all $t>0$, $r(P^D_t)=\mathrm{e}^{-\lambda_0t}$. 

\medskip \noindent \textbf{Step 3}. Let us prove that $\lambda_0>0$. Assume that $\lambda_0=0$, then for all $(q,p)\in D$, $t>0$,
\begin{align*}
    P^D_t\phi(q,p)=\mathbb{E}_{(q,p)}\left[\mathbb{1}_{\tau_\partial>t}\phi(X_t)\right]=\phi(q,p) .
\end{align*}
Hence, $\sup_{q\in\mathcal{O}}\phi(q,p)\leq\Vert\phi\Vert_\infty \sup_{q\in\mathcal{O}}\mathbb{P}_{(q,p)}(\tau_\partial>t)\underset{\vert p\vert\rightarrow\infty}{\longrightarrow}0$ by Lemma~\ref{conv vitesse infinie}. As a consequence, $\phi\in\mathcal{C}^b(\overline{D})$ attains its maximum $\|\phi\|_\infty$ at some $x_0\in\overline{D}$. Then, $\Vert\phi\Vert_\infty=\phi(x_0)\leq\Vert\phi\Vert_\infty \mathbb{P}_{x_0}(\tau_\partial>t)$. Hence, $\mathbb{P}_{x_0}(\tau_\partial>t)=1$, which leads to a contradiction since $\mathbb{P}_{x_0}(X_t\in\mathbb{R}^{2d}\setminus D)>0$. 

\medskip \noindent \textbf{Step 4}. Let us finally prove the properties on $\phi$ stated in Theorem~\ref{properties phi}. First, $\phi\in\mathcal{C}^b(\overline{D})\cap\mathrm{L}^1(D)$ by Proposition~\ref{spectral radius lemma}. In addition, $\phi>0$ on $D\cup\Gamma^-$ and $\phi=0$ on $\Gamma^+\cup\Gamma^0$ using Theorem~\ref{thm density intro qsd} and the fact that $P^D_t\phi=\mathrm{e}^{-\lambda_0t}\phi$. Last, letting $u(t,x)=P^D_t\phi(x)$ in Theorem~\ref{thm density intro qsd}, we get that $u \in \mathcal{C}^\infty(\mathbb{R}_+^* \times D)$ and $\partial_t u = \mathcal{L}u$, but since $u(t,x)$ also writes $\mathrm{e}^{-\lambda_0t}\phi(x)$, we conclude that $\phi\in\mathcal{C}
^\infty(D)$ and $\mathcal{L}\phi=-\lambda_0 \phi$ on $D$.
\end{proof}  
\subsection{Proof of Theorem~\ref{density qsd thm} and Theorem~\ref{vp droite thm}}\label{section qsd spectral interp}
Let us now prove the existence of a unique QSD on the domain $D$ for the processes $(X_t)_{t\geq0}$ in~\eqref{Langevin qsd} and $(\tilde{X}_t)_{t\geq0}$ in~\eqref{Langevin adjoint qsd}.
\begin{proof}[Proof of Theorem~\ref{density qsd thm}]
We prove this theorem for the process $(X_t)_{t\geq0}$ and notice that the exact same proof with the process $(\Tilde{X}_t)_{t\geq0}$ instead of $(X_t)_{t\geq0}$ yields the result for $(\Tilde{X}_t)_{t\geq0}$ with the function $\phi$ instead of $\psi$.

\medskip \noindent\textbf{Step 1}. Let us prove that the measure $\mu$ defined in~\eqref{densite qsd} is a QSD on $D$ of the Langevin process $(X_t)_{t\geq0}$.
 
 Let $t>0$, $A\in\mathcal{B}(D)$. Integrating the equality $\Tilde{P}^D_t\psi=\mathrm{e}^{-(\lambda_0+d\gamma) t}\psi $ over $A$, one obtains that 
 \begin{equation}\label{Step 1 Proof thm density qsd thm}
     \int_A\Tilde{P}^D_t\psi(x)\mathrm{d}x=\mathrm{e}^{-(\lambda_0+d\gamma) t}\int_A\psi(x)\mathrm{d}x.
 \end{equation}
 Furthermore, using Fubini-Tonelli's theorem along with Theorem~\ref{duality thm qsd}, one has that
 $$\int_A\Tilde{P}^D_t\psi(x)\mathrm{d}x=\mathrm{e}^{-d\gamma t}\int_D\psi(x)\mathbb{P}_x(X_t\in A,\tau_\partial>t)\mathrm{d}x.$$
 Therefore, reinjecting into~\eqref{Step 1 Proof thm density qsd thm} we obtain since $\mu(\mathrm{d}x)=\psi(x)\mathrm{d}x$ that
 $$\mathbb{P}_\mu(X_t\in A,\tau_\partial>t)=\mathrm{e}^{-\lambda_0 t}\mu(A),$$
 which gives in particular for $A=D$ that $\mathbb{P}_\mu(\tau_\partial>t)=\mathrm{e}^{-\lambda_0 t}$ and thus $\mu$ is a QSD on $D$ for the process $(X_t)_{t \geq 0}$ by Definition~\ref{def QSD 1}. 
 
\medskip \noindent\textbf{Step 2}. Let $\check{\mu}$ be a QSD on $D$ for the process $(X_t)_{t \geq 0}$. Let us prove that $\check{\mu}=\mu$, where $\mu$ is defined in~\eqref{densite qsd}. We start by proving that $\check{\mu}$ admits a density with respect to the Lebesgue measure on $D$ and that its density is an eigenvector of the semigroup $(\Tilde{P}^D_t)_{t>0}$. By Definition~\ref{def QSD 1}, for all $A\in\mathcal{B}(D)$ and $t>0$,
\begin{equation}\label{qsd uniqueness}
    \mathbb{P}_{\check{\mu}}(X_t\in A, \tau_\partial>t)=\mathbb{P}_{\check{\mu}}(\tau_\partial>t) \check{\mu}(A).
\end{equation}
Moreover, by Proposition~\ref{exit event general case} and the positivity of the transition density $\mathrm{p}^D_t$ on $D\times D$ stated in Theorem~\ref{thm density intro qsd}, there exists $\check{\lambda}_0\in[0,\infty)$ such that $\mathbb{P}_{\check{\mu}}(\tau_\partial>t)=\mathrm{e}^{-\check{\lambda}_0 t}$. 

Let $A\in\mathcal{B}(D)$ with zero Lebesgue measure, then for all $x\in D$, $\mathbb{P}_x(X_t\in A, \tau_\partial>t)=0$ by Theorem~\ref{thm density intro qsd}. As a result, $\mathbb{P}_{\check{\mu}}(X_t\in A, \tau_\partial>t)=0$ and $\check{\mu}(A)=0$ by~\eqref{qsd uniqueness}. Therefore, by Radon-Nikodym's theorem, $\check{\mu}$ admits a measurable non-negative density $\check{\psi}$ with respect to the Lebesgue measure on $D$. Therefore, by~\eqref{qsd uniqueness}, for all $t>0$, $A\in\mathcal{B}(D)$,
\begin{equation}\label{eq1}
    \int_D\check{\psi}(x) \mathbb{P}_x(X_t\in A, \tau_\partial>t) \mathrm{d}x=\mathrm{e}^{-\check{\lambda}_0 t} \int_A \check{\psi}(y) \mathrm{d}y.
\end{equation}
By  Fubini-Tonelli's theorem, 
\begin{align*}
    \int_D\check{\psi}(x) \mathbb{P}_x(X_t\in A, \tau_\partial>t) \mathrm{d}x&=\int_D\int_D\check{\psi}(x) \mathrm{p}^D_t(x,y) \mathbb{1}_A(y) \mathrm{d}y \mathrm{d}x\\
    &=\int_D\mathbb{1}_A(y) \left(\int_D\check{\psi}(x) \mathrm{p}^D_t(x,y) \mathrm{d}x\right) \mathrm{d}y .
\end{align*}
As a result, it follows from~\eqref{eq1} that for almost every $y\in D$, $$\int_D\check{\psi}(x) \mathrm{p}^D_t(x,y) \mathrm{d}x=\check{\psi}(y) \mathrm{e}^{-\check{\lambda}_0 t} .$$
Then, Theorem~\ref{duality thm qsd} ensures that 
$$\int_D\check{\psi}(x) \tilde{\mathrm{p}}^D_t(y,x) \mathrm{d}x=\check{\psi}(y) \mathrm{e}^{-(\check{\lambda}_0+d\gamma) t},$$
which can be rewritten as: for almost every $y\in D$,
\begin{equation}\label{eq5}
     \Tilde{P}^D_t\check{\psi}(y)=\check{\psi}(y) \mathrm{e}^{-(\check{\lambda}_0+d\gamma) t}.
\end{equation} 
Since $\check{\psi}\in\mathrm{L}^1(D)$ and $\tilde{P}^D_t:\mathrm{L}^1(D)\to\mathcal{C}^b(\overline{D})$ by Remark~\ref{semigroupe adjoint}, then $\check{\psi}$ can be chosen in $\mathcal{C}^b(\overline{D})$ by~\eqref{eq5} and in particular~\eqref{eq5} holds for all $y \in \overline{D}$. 

\medskip \noindent \textbf{Step 3}. Let us prove that $\check{\lambda}_0=\lambda_0$. In order to do that let us compute the integral $\int_D\Tilde{P}^D_t\check{\psi}(x)\phi(x)\mathrm{d}x$ in two different ways. On the one hand, it follows from~\eqref{eq5} that
$$\int_D\Tilde{P}^D_t\check{\psi}(x)\phi(x)\mathrm{d}x=\mathrm{e}^{-(\check{\lambda}_0+d\gamma) t}\int_D \check{\psi}(x)\phi(x)\mathrm{d}x.$$
On the other hand, using Fubini-Tonelli's theorem, Theorem~\ref{duality thm qsd} and Theorem~\ref{properties phi}, 
\begin{align*}
\int_D\Tilde{P}^D_t\check{\psi}(x)\phi(x)\mathrm{d}x&=\mathrm{e}^{-d\gamma t}\int_D\check{\psi}(x)P^D_t\phi(x)\mathrm{d}x\\
&=\mathrm{e}^{-(\lambda_0+d\gamma)t}\int_D\check{\psi}(x)\phi(x)\mathrm{d}x.
\end{align*}
Since $\phi>0$, $\check{\psi}\geq0$ on $D$ and $\check{\psi}$ satisfies $\int_D\check{\psi}(x)\mathrm{d}x=1$, then $\int_D\check{\psi}(x)\phi(x)\mathrm{d}x>0$ and $\check{\lambda}_0=\lambda_0$. 

\medskip \noindent \textbf{Step 4}. For $t>0$, $r(\Tilde{P}^D_t)=\mathrm{e}^{-(\lambda_0+d\gamma)t}$ is a simple eigenvalue of $\Tilde{P}^D_t$ (seen as an operator on $\mathcal{C}^b(\overline{D})$) with eigenvector $\psi$ by Theorem~\ref{properties phi}. As a result, since $\int_D\check{\psi}(x)\mathrm{d}x=\int_D\psi(x)\mathrm{d}x=1$, then $\check{\psi}=\psi$ and thus $\check{\mu}=\mu$.
\end{proof}

Let us now prove Theorem~\ref{prop: exit event qsd}.
\begin{proof}[Proof of Theorem~\ref{prop: exit event qsd}]
We shall only prove the first equality, the reasoning is exactly similar for the second equality regarding the QSD of the adjoint process $(\tilde{X}_t)_{t\geq0}$. It is also sufficient to prove the result for $f\in\mathcal{C}^b(\partial D)$, the extension to $f\in\mathrm{L}^\infty(\partial D)$ being straightforward using a dominated convergence argument. Besides, Tietze extension's theorem guarantees that $f\in\mathcal{C}^b(\partial D)$ can be extended to a continuous bounded function on $\mathbb{R}^{2d}$. Such a function can then be approached by a smooth and compactly supported function on $\mathbb{R}^{2d}$. Therefore, it is enough to look at the case $f\in\mathcal{C}_c^\infty(\overline{D})$ and conclude with a dominated convergence argument.

Let $\Phi$ be a non-negative $\mathcal{C}^2$ function in $\mathbb{R}_+$ such that $\Phi(0)=1, \Phi(\rho)=0$ for $\rho\geq1$ and for all $\rho\geq0$, $\Phi(\rho)\in[0,1]$.

Since $\mathcal{O}$ is $\mathcal{C}^2$ there exists $\delta>0$ small enough such that the Euclidean distance to the boundary $\partial\mathcal{O}$, $\mathrm{dist}(\cdot,\partial\mathcal{O})$, is $\mathcal{C}^2$ on the set $\mathcal{O}_\delta:=\{q\in\overline{\mathcal{O}}:\mathrm{dist}(\cdot,\partial\mathcal{O})\leq\delta\}$. Let $\mathrm{d}_\partial$ be defined as the distance function $\mathrm{dist}(\cdot,\partial\mathcal{O})$ on $\mathcal{O}_\delta$ and extended to a $\mathcal{C}^2$ function on $\overline{\mathcal{O}}$. 

For $k\geq1$, we define the function 
\begin{equation}\label{def g_n}
    g_k:(q,p)\in\overline{D}\mapsto f(q,p)\Phi(k\mathrm{d}_\partial(q))\in\mathcal{C}^2_c(\overline{D}).
\end{equation}

For $t>0$, let us compute the limit when $k\rightarrow\infty$ of the following integral in two different ways,
$$\int_D\psi(q,p)\mathbb{E}_{(q,p)}\left[g_k(q_{\tau_\partial\land t},p_{\tau_\partial\land t})\right]\mathrm{d}q\mathrm{d}p.$$

\textbf{Step 1:} For $(q,p)\in D$,
$$\mathbb{E}_{(q,p)}\left[g_k(q_{\tau_\partial\land t},p_{\tau_\partial\land t})\right]=\mathbb{E}_{(q,p)}\left[f(q_{\tau_\partial},p_{\tau_\partial})\mathbb{1}_{\tau_\partial\leq t}\right]+\mathbb{E}_{(q,p)}\left[g_k(q_t,p_t)\mathbb{1}_{\tau_\partial>t}\right].$$
Notice that the sequence $(g_k)_{k\geq1}$ is uniformly bounded by $\Vert f\Vert_\infty$ and that $g_k(q,p)$ converges to $0$ everywhere in $D$. Therefore, by the dominated convergence theorem, since $\psi$ is in $\mathrm{L}^1(D)$,
$$\int_D\psi(q,p)\mathbb{E}_{(q,p)}\left[g_k(q_t,p_t)\mathbb{1}_{\tau_\partial>t}\right]\underset{k\rightarrow\infty}{\longrightarrow}0.$$
As a result,
\begin{equation}\label{lim 1 k}
    \int_D\psi(q,p)\mathbb{E}_{(q,p)}\left[g_k(q_{\tau_\partial\land t},p_{\tau_\partial\land t})\right]\mathrm{d}q\mathrm{d}p\underset{k\rightarrow\infty}{\longrightarrow}\int_D\psi(q,p)\mathbb{E}_{(q,p)}\left[f(q_{\tau_\partial},p_{\tau_\partial})\mathbb{1}_{\tau_\partial\leq t}\right]\mathrm{d}q\mathrm{d}p.
\end{equation}
\textbf{Step 2:} Second, by the Itô formula, we have that
\begin{equation}\label{ito g_n}
    \mathbb{E}_{(q,p)}\left[g_k(q_{\tau_\partial\land t},p_{\tau_\partial\land t})\right]=g_k(q,p)+\mathbb{E}_{(q,p)}\left[\int_0^t\mathcal{L}g_k(q_s,p_s)\mathbb{1}_{\tau_\partial>s}\mathrm{d}s\right].
\end{equation}
By definition of $g_k$~\eqref{def g_n} and $\mathcal{L}$~\eqref{generateur Langevin qsd}, for all $(q,p)\in\mathcal{O}\times\mathbb{R}^d$,
\begin{equation}\label{Lg_n}
    \mathcal{L}g_k(q,p)=k p\cdot\nabla_q d_\partial(q)\Phi'(k\mathrm{d}_\partial(q)) f(q,p) + \Phi(k\mathrm{d}_\partial(q)) \mathcal{L}f(q,p).
\end{equation}
Besides, by definition of $\Phi$,
$$\Phi(k\mathrm{d}_\partial(q))\leq \mathbb{1}_{[0,1/k]}(\mathrm{d}_\partial(q)).$$
As a result, integrating in the first term in the right-hand side of~\eqref{ito g_n} against $\psi$ , 
\begin{equation}\label{int psi g_k}
    \left\vert\int_D\psi(q,p)g_k(q,p)\mathrm{d}q\mathrm{d}p\right\vert\leq\Vert f\Vert_\infty\int_D\psi(q,p)\mathbb{1}_{[0,1/k]}(\mathrm{d}_\partial(q))\mathrm{d}q\mathrm{d}p.
\end{equation}
Furthermore, considering the second term in the right-hand side of the equality~\eqref{Lg_n} one has since $f\in\mathcal{C}_c^\infty(\overline{D})$,
\begin{align*}
    &\left\vert\int_D\psi(q,p)\mathbb{E}_{(q,p)}\left[\int_0^t\Phi(k\mathrm{d}_\partial(q_s)) \mathcal{L}f(q_s,p_s)\mathbb{1}_{\tau_\partial>s}\mathrm{d}s\right]\mathrm{d}q\mathrm{d}p\right\vert\\
    &\leq\Vert\mathcal{L}f\Vert_\infty\int_D\psi(q,p)\mathbb{E}_{(q,p)}\left[\int_0^t\mathbb{1}_{[0,1/k]}(\mathrm{d}_\partial(q_s))\mathbb{1}_{\tau_\partial>s}\mathrm{d}s\right]\mathrm{d}q\mathrm{d}p\\
    &=\Vert\mathcal{L}f\Vert_\infty\int_0^t\int_D\psi(q,p)\mathbb{E}_{(q,p)}\left[\mathbb{1}_{[0,1/k]}(\mathrm{d}_\partial(q_s))\mathbb{1}_{\tau_\partial>s}\right]\mathrm{d}q\mathrm{d}p\mathrm{d}s\\
    &=\Vert\mathcal{L}f\Vert_\infty\frac{1-\mathrm{e}^{-\lambda_0t}}{\lambda_0}\int_D\psi(q,p)\mathbb{1}_{[0,1/k]}(\mathrm{d}_\partial(q))\mathrm{d}q\mathrm{d}p\numberthis\label{int psi ito}  
\end{align*}
since $\psi$ is the QSD on $D$ of $(q_t,p_t)_{t\geq0}$. It remains to estimate the integral against $\psi$ of the first term in the right-hand side of the equality~\eqref{Lg_n}
\begin{align*}
    &\int_D\psi(q,p)\mathbb{E}_{(q,p)}\left[\int_0^tk p\cdot\nabla_q d_\partial(q_s)\Phi'(k\mathrm{d}_\partial(q_s)) f(q_s,p_s)\mathbb{1}_{\tau_\partial>s}\mathrm{d}s\right]\mathrm{d}q\mathrm{d}p\\ 
    &=\frac{1-\mathrm{e}^{-\lambda_0t}}{\lambda_0}\int_D\psi(q,p)k p\cdot\nabla_q d_\partial(q)\Phi'(k\mathrm{d}_\partial(q)) f(q,p)\mathrm{d}q\mathrm{d}p.\numberthis\label{int psi ito 2} 
\end{align*}

Let us study the limit when $k$ goes to infinity of
$$\int_D\psi(q,p)k p\cdot\nabla_q d_\partial(q)\Phi'(k\mathrm{d}_\partial(q)) f(q,p)\mathrm{d}q\mathrm{d}p.$$ 

Recall that $\Phi'$ has support in $[0,1]$. Besides, $\mathrm{d}_\partial$ is $\mathcal{C}^2$ on $\overline{\mathcal{O}}$ and coincides with the Euclidean distance to the boundary on $\mathcal{O}_\delta$. Therefore, for $k\geq1/\delta$, one can apply the following change of variable involving the Weingarten map $W_{\overline{q}}$ and detailed for example in~\cite[Lemma 14.16]{GT},
\begin{align*}
    &\int_D\psi(q,p)k p\cdot\nabla_q d_\partial(q)\Phi'(k\mathrm{d}_\partial(q)) f(q,p)\mathrm{d}q\mathrm{d}p\\
    &=\int_{\partial D}\int_{0}^{1/k}\psi(\overline{q}-\lambda n(\overline{q}),p)k p\cdot\nabla_q d_\partial(\overline{q}-\lambda n(\overline{q}))\Phi'(k\mathrm{d}_\partial(\overline{q}-\lambda n(\overline{q}))) f(\overline{q}-\lambda n(\overline{q}),p)\det(I+\lambda W_{\overline{q}})\mathrm{d}\lambda\sigma_{\partial \mathcal{O}}(\mathrm{d}\overline{q})\mathrm{d}p.\numberthis\label{eq0} 
\end{align*}
Moreover, by the change of variable $s\rightarrow \lambda k$,
\begin{align*}
    &\int_{0}^{1/k}\psi(\overline{q}-\lambda n(\overline{q}),p)k p\cdot\nabla_q d_\partial(\overline{q}-\lambda n(\overline{q}))\Phi'(k\mathrm{d}_\partial(\overline{q}-\lambda n(\overline{q}))) f(\overline{q}-\lambda n(\overline{q}),p)\det(I+\lambda W_{\overline{q}})\mathrm{d}\lambda\\
&=\int_{0}^{1}\psi(\overline{q}-s n(\overline{q})/k,p)p\cdot\nabla_q d_\partial(\overline{q}-s n(\overline{q})/k)\Phi'(k\mathrm{d}_\partial(\overline{q}-s n(\overline{q})/k)) f(\overline{q}-s n(\overline{q})/k,p)\det(I+s W_{\overline{q}}/k)\mathrm{d}s\numberthis\label{eq:weingarten} 
\end{align*}
All the terms inside the integral are bounded. Therefore, using the dominated convergence theorem one easily shows that
\begin{align*}
    &\int_{0}^{1}\psi(\overline{q}-s n(\overline{q})/k,p)p\cdot\nabla_q d_\partial(\overline{q}-s n(\overline{q})/k)\Phi'(k\mathrm{d}_\partial(\overline{q}-s n(\overline{q})/k)) f(\overline{q}-s n(\overline{q})/k,p)\det(I+s W_{\overline{q}}/k)\mathrm{d}s\\
&=\psi(\overline{q},p)f(\overline{q},p)p\cdot \nabla_q\mathrm{d}_\partial(\overline{q})\int_{0}^{1}\Phi'(k\mathrm{d}_\partial(\overline{q}-s n(\overline{q})/k))\mathrm{d}s + \underset{k\rightarrow\infty}{o(k)}\numberthis\label{eq2}. 
\end{align*}
Besides, for $\overline{q}\in\partial\mathcal{O}$, $\nabla_q\mathrm{d}_\partial(\overline{q})=-n(\overline{q})$. In addition, for $s\in(0,1)$,
$$k\mathrm{d}_\partial(\overline{q}-sn(\overline{q})/k)\underset{k\rightarrow\infty}{\longrightarrow}\nabla_q\mathrm{d}_\partial(\overline{q})\cdot (-sn(\overline{q}))=s.$$
As a result,
$$\int_{0}^{1}\Phi'(k\mathrm{d}_\partial(\overline{q}-s n(\overline{q})/k))\mathrm{d}s\underset{k\rightarrow\infty}{\longrightarrow}\int_{0}^{1}\Phi'(s)\mathrm{d}s=\Phi(1)-\Phi(0)=-1.$$
Then, one concludes from~\eqref{eq0},\eqref{eq:weingarten},\eqref{eq2} that 
\begin{equation}\label{int psi bord}
    \int_D\psi(q,p)k p\cdot\nabla_q d_\partial(q)\Phi'(k\mathrm{d}_\partial(q)) f(q,p)\mathrm{d}q\mathrm{d}p\underset{k\rightarrow\infty}{\longrightarrow}\int_{\partial D}\psi(\overline{q},p)f(\overline{q},p)p\cdot n(\overline{q})\sigma_{\partial \mathcal{O}}(\mathrm{d}\overline{q})\mathrm{d}p.
\end{equation}
Consequently, it follows from~\eqref{ito g_n},\eqref{Lg_n},\eqref{int psi g_k},\eqref{int psi ito},\eqref{int psi ito 2} and~\eqref{int psi bord} that
\begin{equation}\label{lim 2 k}
    \int_D\psi(q,p)\mathbb{E}_{(q,p)}\left[g_k(q_{\tau_\partial\land t},p_{\tau_\partial\land t})\right]\mathrm{d}q\mathrm{d}p\underset{k\rightarrow\infty}{\longrightarrow}\frac{1-\mathrm{e}^{-\lambda_0t}}{\lambda_0}\int_{\partial D}\psi(\overline{q},p)f(\overline{q},p)p\cdot n(\overline{q})\sigma_{\partial \mathcal{O}}(\mathrm{d}\overline{q})\mathrm{d}p.
\end{equation}
Finally, it follows from~\eqref{lim 1 k} and~\eqref{lim 2 k} that
$$\int_D\psi(q,p)\mathbb{E}_{(q,p)}\left[f(q_{\tau_\partial},p_{\tau_\partial})\mathbb{1}_{\tau_\partial\leq t}\right]\mathrm{d}q\mathrm{d}p=\frac{1-\mathrm{e}^{-\lambda_0t}}{\lambda_0}\int_{\partial D}\psi(\overline{q},p)f(\overline{q},p)p\cdot n(\overline{q})\sigma_{\partial \mathcal{O}}(\mathrm{d}\overline{q})\mathrm{d}p$$
which concludes the proof since $\psi=0$ on $\Gamma^-$ and positive on $\Gamma^+$ by Theorem~\ref{properties phi}.
\end{proof}

Let us conclude this subsection by proving Theorem~\ref{vp droite thm}. This will provide a spectral interpretation of the QSD on $D$ of the Langevin process, similarly to the spectral interpretation obtained in the overdamped Langevin case, cf. Theorem~\ref{qsd overdamped}.

\begin{proof}[Proof of Theorem~\ref{vp droite thm}]
The couple $(\lambda_0,\phi)$, defined in Theorem~\ref{properties phi}, is clearly a solution to the left eigenvalue problem in~\eqref{vp droite gauche}. Let us prove that such a couple $(\lambda_0,\phi)$ is unique, up to a multiplicative constant for $\phi$. Since the reasoning for the right eigenvalue problem with solution $(\lambda_0,\psi)$ is the same, with the process $(\tilde{X}_t)_{t\geq0}$ instead of $(X_t)_{t\geq0}$, it will not be detailed. 

Let $\lambda\in\mathbb{R}$ and  $\eta\in\mathcal{C}^{2}(D)\cap\mathcal{C}^{b}(D\cup\Gamma^+)$ be a non-zero and non-negative classical solution of the left eigenvalue problem in~\eqref{vp droite gauche}. Let $$\tau_{V^c_k}:=\inf\{t>0: X_t\notin V_k\},$$ where $V_k:=\{(q,p)\in D : \vert p\vert< k, \mathrm{d}_\partial(q)>\frac{1}{k} \}$ and we recall that $\mathrm{d}_\partial$ refers to the distance to $\partial \mathcal{O}$. Applying Itô's formula to the process $(\mathrm{e}^{\lambda s} \eta(X_s))_{s\geq0}$ at the stopping time $t\land\tau_{V^c_k}$, one gets, for $x\in D$, $\mathbb{P}_x$-almost surely, for all $t\geq0$, 
\begin{equation}\label{ito formula}
    \mathrm{e}^{\lambda (t\land\tau_{V^c_k})} \eta\big(X_{t\land\tau_{V^c_k}}\big)=\eta(x)+\sigma \int_0^t\mathbb{1}_{s\leq\tau_{V^c_k}} \mathrm{e}^{\lambda s} \nabla_p\eta(X_s)\cdot \mathrm{d}B_s,
\end{equation}
since $\mathcal{L}\eta+\lambda \eta=0$ on $D$. Moreover, $\nabla_p\eta$ is bounded on the compact $\overline{V_k}$ since $\eta\in\mathcal{C}^{2}(D)$. Therefore, the stochastic integral in the right-hand side of the equality~\eqref{ito formula} is a martingale and its expectation vanishes. Hence,
\[\begin{aligned}
\mathbb{E}_x\bigg[\mathrm{e}^{\lambda (t\land\tau_{V^c_k})} \eta\big(X_{t\land\tau_{V^c_k}}\big)\bigg]&=\eta(x),
\end{aligned}\]
which can be rewritten as
\begin{equation}\label{eq3}
    \eta(x)=\mathrm{e}^{\lambda t} \mathbb{E}_x\bigg[ \eta\big(X_{t}\big) \mathbb{1}_{\tau_{V^c_k}>t} \bigg]+\mathbb{E}_x\bigg[\mathrm{e}^{\lambda \tau_{V^c_k}} \eta\big(X_{\tau_{V^c_k}}\big) \mathbb{1}_{\tau_{V^c_k}\leq t}\bigg].
\end{equation}

Now we would like to let $k\rightarrow\infty$. Let us prove the following limit, $\mathbb{P}_x$-almost surely,
$$\lim_{k\rightarrow\infty}\tau_{V_k^c}=\tau_{\partial},$$
using the same reasoning as in~\cite[Section 3.1]{LelRamRey}.

The sequence $(\tau_{V_k^c})_{k\geq1}$ is an increasing sequence of random variables, therefore it converges almost surely to $\sup_{k\geq1}\tau_{V_k^c}$. Besides, using the continuity of the trajectories of $(X_t)_{t\geq0}$, one gets for all $r>0$,
\[\begin{aligned}
\left\{\sup_{k\geq1}\tau_{V_k^c}>r\right\}&=\left\{\exists  k\geq1: \tau_{V_k^c}>r\right\}\\
&=\left\{\exists  k\geq1: \sup_{u\in[0,r]}\vert p_u\vert< k, \inf_{u\in[0,r]} \mathrm{d}_\partial(q_u)>\frac{1}{k}\right\}\\
&=\left\{\sup_{u\in[0,r]}\vert p_u\vert<\infty, \inf_{u\in[0,r]} \mathrm{d}_\partial(q_u)>0\right\}\\
&=\left\{\sup_{u\in[0,r]}\vert p_u\vert<\infty, \tau_\partial>r\right\} .
\end{aligned}\]
 
For all $r>0$, we have that $\sup_{u\in[0,r]}\vert p_u\vert<\infty$, almost surely. Therefore, $\sup_{k\geq1}\tau_{V_k^c}>r$ if and only if $\tau_\partial>r$, that is to say $\sup_{k\geq1}\tau_{V_k^c}=\tau_\partial$ almost surely.
As a result, one gets $\lim_{k\rightarrow\infty}\tau_{V_k^c}=\tau_{\partial}$ almost surely. Besides, since $(\tau_{V^c_k})_{k\geq1}$ is an increasing sequence, then for all $s>0$, almost surely,
$$\mathbb{1}_{\tau_{V_k^c}> s}\underset{k\rightarrow\infty}{\longrightarrow} \mathbb{1}_{\tau_{\partial}> s}.$$
As a result, by continuity of the trajectories of $(X_t)_{t\geq0}$, $\mathbb{P}_x$-almost surely, 
$$\mathbb{1}_{\tau_{V_k^c}\leq t} \eta(X_{\tau_{V_k^c}})\underset{k\rightarrow\infty}{\longrightarrow} \mathbb{1}_{\tau_{\partial}\leq t} \eta(X_{\tau_{\partial}}).$$
Moreover $\eta(X_{\tau_{\partial}})=0$ almost surely on the event $\{\tau_{\partial}\leq t\}$ since $X_{\tau_{\partial}}\in\Gamma^+$ $\mathbb{P}_x$-almost surely by~\cite[Assertion (ii) of Proposition 2.8]{LelRamRey}. Therefore, taking the limit $k\rightarrow\infty$ in~\eqref{eq3}, one gets by the dominated convergence theorem that
\begin{equation}\label{eq4}
    \forall t>0,\quad\forall x\in D,\qquad\mathbb{E}_x\big[ \eta\big(X_{t}\big) \mathbb{1}_{\tau_{\partial}>t} \big]=\mathrm{e}^{-\lambda t} \eta(x), 
\end{equation}
which ensures in particular that necessarily $\lambda\geq0$. This also writes 
$$\forall t>0,\quad\forall x\in D,\qquad\int_D\mathrm{p}^D_t(x,y)\eta(y)\mathrm{d}y=\mathrm{e}^{-\lambda t} \eta(x) . $$
Using the boundedness of $\eta$ along with~\eqref{it:compacite-qsd:2} in Theorem~\ref{compactness}, we deduce that $\eta\in\mathrm{L}^1(D)$. Now let $\tilde{\eta}=\eta/\int_D\eta$, then using Theorem~\ref{duality thm qsd} one obtains that
$$\forall t>0,\quad\forall x\in D,\qquad\int_D\tilde{\eta}(y)\tilde{\mathrm{p}}^D_t(y,x)\mathrm{d}y=\mathrm{e}^{-(\lambda+d\gamma) t} \tilde{\eta}(x) . $$
Integrating over $D$ we obtain that $\mathbb{P}_{\tilde{\nu}}(\tau_{\partial}>t)=\mathrm{e}^{-(\lambda+d\gamma) t}$ with $\tilde{\nu}(\mathrm{d}x)=\tilde{\eta}(x)\mathrm{d}x$. Then, integrating  over any $A\in\mathcal{B}(D)$, we obtain that 
$\tilde{\nu}$ is a QSD on $D$ of the process $(\tilde{X}_t)_{t\geq0}$. Consequently, the uniqueness of such a QSD, by Theorem~\ref{density qsd thm}, ensures that $\tilde{\nu}=\tilde{\mu}$ where $\tilde{\mu}$ is defined in Theorem~\ref{density qsd thm}. In addition, it implies that $\lambda=\lambda_0$, which concludes the proof.  
\end{proof}

\subsection{Proof of Theorems~\ref{spectral decomp} and~\ref{cv semigroupe conditionnel}}\label{section long time semigroupe}

This section is devoted to the study of the long-time convergence of the semigroup $(P^D_t)_{t>0}$. Note that a similar study could be performed for the semigroup $(\Tilde{P}^D_t)_{t>0}$, using the duality between the two semigroups as stated in Theorem~\ref{duality thm qsd}.  

We start this subsection by ensuring the existence of a spectral gap for the operator $P^D_t$. In the next statement, we denote by $\mathcal{C}^b(\overline{D},\mathbb{C})$ the space of complex-valued continuous bounded functions on $\overline{D}$.
\begin{lemma}[Spectral gap]\label{Lemma radius eigenvalue} Under Assumptions~\ref{hyp F1 qsd} and~\ref{hyp O qsd}, for all $t>0$, the operator $P^D_t$ admits a unique complex eigenvalue with modulus equal to $r(P^D_t)=\mathrm{e}^{-\lambda_0 t}$ and eigenvector in $\mathcal{C}^b(\overline{D},\mathbb{C})$. 
\end{lemma}
\begin{proof}
Assume that there exists an eigenvector $h_t\in\mathcal{C}^b(\overline{D},\mathbb{C})$ of $P^D_t$ with eigenvalue $z\in\mathbb{C}\setminus\{\mathrm{e}^{-\lambda_0 t}\}$ such that $\vert z\vert=\mathrm{e}^{-\lambda_0 t}$. Let $\psi\in\mathcal{C}^b(\overline{D})$ be the eigenvector of $\tilde{P}^D_t$ from Theorem~\ref{properties phi}. 

First, let us prove that $\int_D h_t(x)\psi(x)\mathrm{d}x=0$ by computing the integral $\int_D h_t(x)\Tilde{P}^D_t\psi(x)\mathrm{d}x$ in two different ways. On the one hand, by Theorem~\ref{properties phi}, 
$$\int_D h_t(x)\Tilde{P}^D_t\psi(x)\mathrm{d}x=\mathrm{e}^{-(\lambda_0+d\gamma) t}\int_D h_t(x)\psi(x)\mathrm{d}x  .$$
On the other hand, by Theorem~\ref{duality thm qsd},  
\begin{align*}
    \int_D h_t(x)\Tilde{P}^D_t\psi(x)\mathrm{d}x&=\mathrm{e}^{-d\gamma t}\int_D P^D_th_t(x)\psi(x)\mathrm{d}x\\
    &=z \mathrm{e}^{-d\gamma t}\int_D h_t(x)\psi(x)\mathrm{d}x .
\end{align*}
Therefore, since $z\neq\mathrm{e}^{-\lambda_0 t}$, $\int_D h_t(x)\psi(x)\mathrm{d}x=0$, and in particular
\begin{equation}\label{eq eigenvector}
    \int_D \mathrm{Re}(h_t(x))\psi(x)\mathrm{d}x=\int_D \mathrm{Im}(h_t(x))\psi(x)\mathrm{d}x=0.
\end{equation}

Besides, one has for $x\in D$,
\begin{align*}
    r(P^D_t)\vert h_t(x)\vert&=\vert P^D_th_t(x)\vert\\
    &=\Big\vert \int_D\mathrm{p}^D_t(x,y)h_t(y)\mathrm{d}y\Big\vert\\
    &< \int_D\mathrm{p}^D_t(x,y)\vert h_t(y)\vert\mathrm{d}y=P^D_t\vert h_t\vert(x) ,
\end{align*}
by the triangle inequality since the equality case requires that $\mathrm{Re}(h_t)$ and $\mathrm{Im}(h_t)$ have constant signs on $D$, which would imply $h_t=0$ from~\eqref{eq eigenvector} since $\psi>0$ on $D$. As a result,
\begin{align*}
    r(P^D_t)\int_D\psi(x) \vert h_t(x)\vert\mathrm{d}x&<\int_D\psi(x) P^D_t\vert h_t\vert(x)\mathrm{d}x\\
    &=\mathrm{e}^{d\gamma t}\int_D\Tilde{P}^D_t\psi(x) \vert h_t(x)\vert\mathrm{d}x\\
    &=r(P^D_t)\int_D\psi(x) \vert h_t(x)\vert\mathrm{d}x 
\end{align*}
which leads to a contradiction, therefore such an eigenvalue does not exist. 
\end{proof} 
We are now able to prove Theorem~\ref{spectral decomp}.

\begin{proof}[Proof of Theorem~\ref{spectral decomp}]
\textbf{Step 1:} We shall first prove~\eqref{long time cv P^D_t} for $f\in\mathcal{C}^b(\overline{D})$ then use a regularization argument of the semigroup $(P^D_t)_{t\geq0}$ in the next step to conclude. 

Let us define the following vector space of $\mathcal{C}^b(\overline{D})$,
$$\mathrm{Span}(\psi)^\perp:=\left\{f\in\mathcal{C}^b(\overline{D}):\int_Df(x)\psi(x)\mathrm{d}x=0\right\}.$$
On the one hand, it is clear that this is a closed subset of $\mathcal{C}^b(\overline{D})$, and thus a Banach space. On the other hand, it follows from Theorems~\ref{duality thm qsd} and~\ref{properties phi} that $\mathrm{Span}(\psi)^\perp$ is stable by $P^D_1$. As a consequence, we may consider in the sequel the operator $P^D_1|_{\mathrm{Span}(\psi)^\perp}$.

The compactness of $P^D_1$ ensures the compactness of  $P^D_1|_{\mathrm{Span}(\psi)^\perp}$ as well. Therefore, any non-zero element of the spectrum $\sigma(P^D_1|_{\mathrm{Span}(\psi)^\perp})$ is an eigenvalue of $P^D_1|_{\mathrm{Span}(\psi)^\perp}$ and the eigenvalues can only accumulate at $0$. Therefore, if the spectral radius $r(P^D_1|_{\mathrm{Span}(\psi)^\perp})>0$, then it is an eigenvalue of $P^D_1|_{\mathrm{Span}(\psi)^\perp}$ by Definition~\ref{def spectral radius}. Moreover, Lemma~\ref{Lemma radius eigenvalue} ensures that $r(P^D_1|_{\mathrm{Span}(\psi)^\perp})<r(P^D_1)$ since $r(P^D_1)$ is a simple eigenvalue associated to a positive function $\phi$ which thus does not belong to $\mathrm{Span}(\psi)^\perp$. 

In any case, we thus have $r(P^D_1|_{\mathrm{Span}(\psi)^\perp})<r(P^D_1)$, so that there exists $\alpha^* \in (0,+\infty]$ such that $r(P^D_1|_{\mathrm{Span}(\psi)^\perp})=\mathrm{e}^{-\lambda_0-\alpha^*}$. In addition, for $\alpha\in[0,\alpha^*)$, by Proposition~\ref{spectral radius} there exists $N_0\geq1$ such that for all $N\geq N_0$,
$$\vertiii{P^D_N|_{\mathrm{Span}(\psi)^\perp}}_{\mathcal{C}^b(\overline{D})}\leq\mathrm{e}^{-(\lambda_0+\alpha)N},$$ where we have used the semigroup property to write $(P^D_1|_{\mathrm{Span}(\psi)^\perp})^N={P^D_N}|_{\mathrm{Span}(\psi)^\perp}$.
 
Noticing that for all $f\in\mathcal{C}^b(\overline{D})$, $f-\frac{\phi\otimes \psi}{\int_D \phi \psi}(f)\in\mathrm{Span}(\psi)^\perp$ since $\phi\otimes \psi(f)=(\int_D \psi f)\phi$, one gets for $N\geq N_0$ and $f\in\mathcal{C}^b(\overline{D})$, 
$$\left\Vert P^D_Nf-\mathrm{e}^{-\lambda_0N}\frac{\phi\otimes \psi}{\int_D \phi \psi}(f)\right\Vert_\infty=\left\Vert P^D_N\left(f-\frac{\phi\otimes \psi}{\int_D \phi \psi}(f)\right)\right\Vert_\infty\leq\mathrm{e}^{-(\lambda_0+\alpha)N} \left\Vert f-\frac{\phi\otimes \psi}{\int_D \phi \psi}(f)\right\Vert_\infty .$$
Let $t\geq N_0$, then $\lfloor t\rfloor\geq N_0$, and we have that
\begin{align*}
    \left\Vert P^D_tf-\mathrm{e}^{-\lambda_0t}\frac{\phi\otimes \psi}{\int_D \phi \psi}(f)\right\Vert_\infty&=\left\Vert P^D_{t-\lfloor t\rfloor}P^D_{\lfloor t\rfloor}\left(f-\frac{\phi\otimes \psi}{\int_D \phi \psi}(f)\right)\right\Vert_\infty\\
    &\leq\left\Vert P^D_{\lfloor t\rfloor}\left(f-\frac{\phi\otimes \psi}{\int_D \phi \psi}(f)\right)\right\Vert_\infty\\
    &\leq\mathrm{e}^{-(\lambda_0+\alpha)\lfloor t\rfloor}\left\Vert f-\frac{\phi\otimes \psi}{\int_D \phi \psi}(f)\right\Vert_\infty\\
    &\leq\mathrm{e}^{-(\lambda_0+\alpha)t} \mathrm{e}^{\lambda_0+\alpha} \left\Vert f-\frac{\phi\otimes \psi}{\int_D \phi \psi}(f)\right\Vert_\infty\\
    &\leq\mathrm{e}^{-(\lambda_0+\alpha)t} \mathrm{e}^{\lambda_0+\alpha}\left(1+\frac{\Vert \phi\Vert_\infty}{\int_D\phi\psi}\right) \left\Vert f\right\Vert_\infty,
\end{align*}
which concludes the proof of~\eqref{long time cv P^D_t} when $f\in\mathcal{C}^b(\overline{D})$, since the behavior of the left-hand side for $t \leq N_0$ can easily be bounded appropriately. 

\noindent\textbf{Step 2:} Let us extend~\eqref{long time cv P^D_t} for $f\in\mathrm{L}^\infty(D)$. It follows from Theorem~\ref{compactness} that for any $s>0$ and any $f\in\mathrm{L}^\infty(D)$, $P^D_s f\in\mathcal{C}^b(\overline{D})$. Therefore, by \textbf{Step 1}, for all $\alpha\in[0,\alpha^*)$, there exists $C_\alpha>0$  such that for all $t>0$, $\epsilon\in(0,t)$ and $f\in\mathrm{L}^\infty(D)$,  
$$\left\Vert P^D_{t-\epsilon}(P^D_\epsilon f)-\mathrm{e}^{-\lambda_0 (t-\epsilon)}\frac{\phi\otimes \psi(P^D_\epsilon f)}{\int_D \phi \psi}\right\Vert_\infty\leq C_\alpha\Vert P^D_\epsilon f\Vert_\infty\mathrm{e}^{-(\lambda_0+\alpha) (t-\epsilon)}.$$
Furthermore, $ P^D_{t-\epsilon}(P^D_\epsilon f)= P^D_t f$, $\Vert P^D_\epsilon f\Vert_\infty\leq\Vert f\Vert_{\mathrm{L}^\infty(D)}$ and by the Fubini permutation along with Theorem~\ref{duality thm qsd} since $f,\psi\in\mathrm{L}^\infty(D)$,  
$$\mathrm{e}^{\lambda_0\epsilon}\phi\otimes \psi(P^D_\epsilon f)=\mathrm{e}^{\lambda_0\epsilon}\phi \int_D\psi P^D_\epsilon f=\mathrm{e}^{(\lambda_0+d\gamma)\epsilon}\phi \int_Df\tilde{P}^D_\epsilon\psi =\phi\otimes \psi(f).$$

As a result, for all $\alpha\in[0,\alpha^*)$, there exists $C_\alpha>0$  such that for all $t>0$, $\epsilon\in(0,t)$ and $f\in\mathrm{L}^\infty(D)$,
$$\left\Vert P^D_t f-\mathrm{e}^{-\lambda_0 t}\frac{\phi\otimes \psi(f)}{\int_D \phi \psi}\right\Vert_\infty\leq C_\alpha\Vert f\Vert_{\mathrm{L}^\infty(D)}\mathrm{e}^{-(\lambda_0+\alpha) (t-\epsilon)},$$
which yields the proof of \textbf{Step 2}.
\end{proof}

We conclude this section with a proof of the long-time convergence, in total variation, of the distribution of the Langevin process conditioned to remain in $D$, towards its QSD on $D$.

\begin{proof}[Proof of Theorem~\ref{cv semigroupe conditionnel}]
We first recall that for any probability measure $\theta$ on $D$ and $t \geq 0$,
\begin{equation*}
  \left\|\mathbb{P}_\theta\left(X_t \in \cdot | \tau_\partial > t\right) - \mu\right\|_{TV} = \sup_{f \in \mathrm{L}^\infty(D), \|f\|_{\mathrm{L}^\infty(D)} \leq 1} \left|\mathbb{E}_\theta\left[f(X_t)|\tau_\partial > t\right] - \int_D f \mathrm{d}\mu\right|.
\end{equation*}
Let us fix $\alpha \in [0,\alpha^*)$ and show that there exists $C'_\alpha$ such that for any initial distribution $\theta$, $t \geq 0$ and $f \in \mathrm{L}^\infty(D)$,
\begin{equation*}
  \left\vert\frac{\mathbb{E}_{\theta}(f(X_t)\mathbb{1}_{\tau_\partial>t})}{\mathbb{P}_{\theta}(\tau_\partial>t)}-\int_D f\psi\right\vert\leq \frac{C'_\alpha}{\int_D \phi \mathrm{d}\theta} \mathrm{e}^{-\alpha t}\Vert f\Vert_{\mathrm{L}^\infty(D)}.
\end{equation*} 

First, let us prove that $$\mathbb{P}_{\theta}(\tau_\partial>t)\geq \frac{\int_D\phi \mathrm{d}{\theta}}{\Vert\phi\Vert_\infty}\mathrm{e}^{-\lambda_0 t}.$$ For $x\in D$, $t>0$ one has
\begin{align*}
    \mathbb{P}_x(\tau_\partial>t)&=\int_D\mathrm{p}^D_t(x,y) \mathrm{d}y\\
    &=\int_D\mathrm{p}^D_t(x,y)\frac{\phi(y)}{\phi(y)}\mathrm{d}y\\
    &\geq \frac{\int_D\mathrm{p}^D_t(x,y)\phi(y)\mathrm{d}y}{\Vert\phi\Vert_\infty}=\frac{\phi(x) \mathrm{e}^{-\lambda_0 t}}{\Vert\phi\Vert_\infty},
\end{align*}
by Theorem~\ref{properties phi}. Therefore,
\begin{equation}\label{minoration proba}
    \mathbb{P}_{\theta}(\tau_\partial>t)=\int_D\theta(\mathrm{d}x)\mathbb{P}_x(\tau_\partial>t)\geq\frac{\int_D\phi \mathrm{d}{\theta}}{\Vert\phi\Vert_\infty}\mathrm{e}^{-\lambda_0 t}. 
\end{equation}
As a result, it follows from Theorem~\ref{spectral decomp} and the inequality~\eqref{minoration proba}, the existence of $C_{\alpha}>0$ such that for all $t>0$,
\begin{align*}
  \left\vert\frac{\mathbb{E}_{\theta}(f(X_t)\mathbb{1}_{\tau_\partial>t})}{\mathbb{P}_{\theta}(\tau_\partial>t)}-\int_D f\psi\right\vert&=\left\vert\frac{\int_D\left(\mathbb{E}_x(f(X_t)\mathbb{1}_{\tau_\partial>t})-\left(\int_D f\psi\right) \mathbb{P}_{\theta}(\tau_\partial>t)\right){\theta}(\mathrm{d}x)}{\mathbb{P}_{\theta}(\tau_\partial>t)}\right\vert\\
  &\leq\int_D\frac{\left\vert\mathbb{E}_x(f(X_t)\mathbb{1}_{\tau_\partial>t})-\mathrm{e}^{-\lambda_0 t}\frac{\int_D \psi f}{\int_D \phi \psi} \phi(x)\right\vert}{\mathbb{P}_{\theta}(\tau_\partial>t)}{\theta}(\mathrm{d}x)\\
  &+\left\vert\int_D \psi f\right\vert\frac{\left\vert\mathrm{e}^{-\lambda_0 t}\frac{\int_D\phi \mathrm{d}{\theta}}{\int_D \phi \psi} -  \mathbb{P}_{\theta}(\tau_\partial>t)\right\vert}{\mathbb{P}_{\theta}(\tau_\partial>t)}\\
  &\leq\frac{C_{\alpha}}{\mathbb{P}_{\theta}(\tau_\partial>t)} \mathrm{e}^{-(\lambda_0+\alpha )t}\Vert f\Vert_{\mathrm{L}^\infty(D)} +\Vert f\Vert_{\mathrm{L}^\infty(D)}\frac{C_{\alpha}\mathrm{e}^{-(\lambda_0+\alpha) t} }{\mathbb{P}_{\theta}(\tau_\partial>t)}\\
  &\leq\frac{2 C_{\alpha}}{\int_D \phi \mathrm{d}{\theta}} \mathrm{e}^{-\alpha t}\Vert f\Vert_{\mathrm{L}^\infty(D)}\Vert \phi\Vert_\infty. \qedhere
\end{align*}
\end{proof}
 
\noindent \textbf{Acknowledgments :} Mouad Ramil is supported by the Région Ile-de- France through a PhD fellowship of the Domaine d'Intérêt Majeur (DIM) Math Innov.  This work also benefited from the support of the projects ANR EFI (ANR-17-CE40-0030) and ANR QuAMProcs (ANR-19-CE40-0010) from the French National Research Agency. Finally, Tony Lelièvre has received funding from the European Research-Council (ERC) under the European Union’s Horizon 2020 research and innovation programme (grant agreement No 810367), project EMC2. The authors would also like to thank Etienne Bernard (CERMICS) for useful discussions on analytical properties of semigroups. 

\bibliographystyle{plain}
\bibliography{biblio}

\begin{thebibliography}{10}

\bibitem{V3}
N.~Champagnat, K.~A. Coulibaly-Pasquier, and D.~Villemonais.
\newblock Criteria for exponential convergence to quasi-stationary
  distributions and applications to multi-dimensional diffusions.
\newblock In {\em S\'{e}minaire de {P}robabilit\'{e}s {XLIX}}, volume 2215 of
  {\em Lecture Notes in Math.}, pages 165--182. Springer, Cham, 2018.

\bibitem{V}
N.~{Champagnat} and D.~{Villemonais}.
\newblock {General criteria for the study of quasi-stationarity}.
\newblock {\em arXiv e-prints}, page arXiv:1712.08092, Dec 2017.

\bibitem{V2}
N.~{Champagnat} and D.~{Villemonais}.
\newblock {Lyapunov criteria for uniform convergence of conditional
  distributions of absorbed Markov processes}.
\newblock {\em arXiv e-prints}, page arXiv:1704.01928, Apr 2017.

\bibitem{Collet}
P.~Collet, S.~Mart\'{\i}nez, and J.~San~Mart\'{\i}n.
\newblock {\em Quasi-stationary distributions}.
\newblock Probability and its Applications (New York). Springer, Heidelberg,
  2013.
\newblock Markov chains, diffusions and dynamical systems.

\bibitem{F}
A.~Friedman.
\newblock {\em Stochastic differential equations and applications. {V}ol. 1}.
\newblock Academic Press, New York-London, 1975.
\newblock Probability and Mathematical Statistics, Vol. 28.

\bibitem{GT}
D.~Gilbarg and N.~S. Trudinger.
\newblock {\em Elliptic partial differential equations of second order}.
\newblock Classics in Mathematics. Springer-Verlag, Berlin, 2001.
\newblock Reprint of the 1998 edition.

\bibitem{GQZ}
G.~L. Gong, M.~P. Qian, and Z.~X. Zhao.
\newblock Killed diffusions and their conditioning.
\newblock {\em Probab. Theory Related Fields}, 80:151--167, 1988.

\bibitem{GrutWid}
M.~Gr{\"u}ter and K.-O. Widman.
\newblock The {G}reen function for uniformly elliptic equations.
\newblock {\em Manuscripta Mathematica}, 37(3):303--342, 1982.

\bibitem{GuiNectoux}
A.~Guillin, B.~Nectoux, and L.~Wu.
\newblock {Quasi-stationary distribution for strongly Feller Markov processes
  by Lyapunov functions and applications to hypoelliptic Hamiltonian systems}.
\newblock {\em https://hal.archives-ouvertes.fr/hal-03068461/}, 2020.

\bibitem{PositiveDensity}
D.~P. Herzog and J.~C. Mattingly.
\newblock A practical criterion for positivity of transition densities.
\newblock {\em Nonlinearity}, 28(8):2823--2845, 2015.

\bibitem{KnobPart}
R.~Knobloch and L.~Partzsch.
\newblock Uniform conditional ergodicity and intrinsic ultracontractivity.
\newblock {\em Potential Anal.}, 33(2):107--136, 2010.

\bibitem{LebLelPer}
C.~Le~Bris, T.~Leli\`evre, M.~Luskin, and D.~Perez.
\newblock A mathematical formalization of the parallel replica dynamics.
\newblock {\em Monte Carlo Methods Appl.}, 18(2):119--146, 2012.

\bibitem{LelRouSto10}
T.~Leli\`evre, M.~Rousset, and G.~Stoltz.
\newblock {\em Free energy computations}.
\newblock Imperial College Press, London, 2010.
\newblock A mathematical perspective.

\bibitem{LelSto16}
T.~Leli\`evre and G.~Stoltz.
\newblock Partial differential equations and stochastic methods in molecular
  dynamics.
\newblock {\em Acta Numer.}, 25:681--880, 2016.

\bibitem{LelRamRey}
T.~Lelièvre, M.~Ramil, and J.~Reygner.
\newblock {A probabilistic study of the kinetic {F}okker-{P}lanck equation in
  cylindrical domains}.
\newblock {\em arXiv e-prints}, page arXiv:2010.10157, Jan 2021.

\bibitem{VilMel}
S.~M\'{e}l\'{e}ard and D.~Villemonais.
\newblock Quasi-stationary distributions and population processes.
\newblock {\em Probab. Surv.}, 9:340--410, 2012.

\bibitem{nier-18}
F.~Nier.
\newblock {\em Boundary conditions and subelliptic estimates for geometric
  Kramers-Fokker-Planck operators on manifolds with boundaries}, volume 252.
\newblock American Mathematical Society, 2018.

\bibitem{perez-uberuaga-voter-15}
D.~Perez, B.P. Uberuaga, and A.F. Voter.
\newblock The parallel replica dynamics method--{C}oming of age.
\newblock {\em Computational Materials Science}, 100:90--103, 2015.

\bibitem{RamPHD}
M.~Ramil.
\newblock {\em Processus cinétiques dans des domaines à bord et
  quasi-stationnarité}.
\newblock PhD thesis, Ecole des Ponts ParisTech, 2020.

\bibitem{RS}
M.~Reed and B.~Simon.
\newblock {\em Methods of modern mathematical physics. {I}. {F}unctional
  analysis}.
\newblock Academic Press, New York-London, 1972.

\bibitem{RB}
L.~Rey-Bellet.
\newblock Ergodic properties of {M}arkov processes.
\newblock In {\em Open quantum systems. {II}}, volume 1881 of {\em Lecture
  Notes in Math.}, pages 1--39. Springer, Berlin, 2006.

\bibitem{SW}
H.H. Schaefer and M.P. Wolff.
\newblock Topological vector spaces.
\newblock {\em Graduate Texts in Mathematics}, 1966.

\bibitem{MSM}
C.~Sch{\"u}tte and M.~Sarich.
\newblock {\em Metastability and Markov State Models in Molecular Dynamics},
  volume~24.
\newblock American Mathematical Soc., 2013.

\bibitem{kMC}
A.~Voter.
\newblock Introduction to the kinetic {M}onte {C}arlo method.
\newblock In {\em Radiation effects in solids}, pages 1--23. Springer, 2007.

\end{thebibliography}

\end{document}